\documentclass[notitlepage,twoside,a4paper]{amsart}
\usepackage{mathrsfs}
\usepackage{amsmath,amssymb,enumerate}
\usepackage{epsfig,fancyhdr,color}
\usepackage{epstopdf}
\usepackage{amssymb}
\usepackage{amsmath,amsthm}
\usepackage{latexsym}
\usepackage{amscd}
\usepackage{psfrag}
\usepackage{graphicx}
\usepackage{epsf}
\usepackage[latin1]{inputenc}
\usepackage[all]{xy}
\usepackage{tikz}
\usepackage{prettyref}
\usepackage{subfigure}
\usepackage{float}
\usepackage{color}
\usepackage{array}
\usepackage{cite}
\usepackage[totalwidth=15cm,totalheight=24cm]{geometry}

\newcommand{\tq}{\, \big| \, }
\newcommand{\II}{I\hspace{-0.1cm}I}
\newcommand{\III}{I\hspace{-0.1cm}I\hspace{-0.1cm}I}

\DeclareMathOperator{\area}{Area}

\newtheorem{theorem}{\rm\bf Theorem}[section]
\newtheorem{proposition}[theorem]{\rm\bf Proposition}

\newtheorem{lemma}[theorem]{\rm\bf Lemma}
\newtheorem{corollary}[theorem]{\rm\bf Corollary}
\newtheorem{definition}[theorem]{\rm\bf Definition}
\newtheorem{remark}[theorem]{\rm\bf Remark}

\newtheorem{question}[theorem]{\rm\bf Question}

\newtheoremstyle{named}{}{}{\itshape}{}{\bfseries}{.}{.5em}{#1 \thmnote{#3}}
\theoremstyle{named}

\newcommand{\C}{{\mathbb C}}
\newcommand{\CP}{{\mathbb CP}}
\newcommand{\N}{{\mathbb N}}

\newcommand{\HH}{{\mathbb H}}
\newcommand{\R}{{\mathbb R}}

\newcommand{\Z}{{\mathbb Z}}
\newcommand{\bCP}{{\mathbb{CP}}}

\newcommand{\cE}{{\mathcal E}}
\newcommand{\cG}{{\mathcal G}}
\newcommand{\cL}{{\mathcal L}}
\newcommand{\cM}{{\mathcal M}}

\newcommand{\cP}{{\mathcal P}}
\newcommand{\cT}{{\mathcal T}}
\newcommand{\cML}{{\mathcal{ML}}}

\newcommand{\Isom}{\rm{Isom}}
\newcommand{\PSL}{\rm{PSL}}

\newcommand{\tr}{{\rm tr}}

\newcommand{\dv}{{\rm div}}

\newcommand{\cc}{\frak{c}}
\newcommand{\cm}{\frak{m}}
\newcommand{\ct}{\frak{t}}

\DeclareMathOperator{\adj}{adj}

\newcounter{notes}%

\def\interieur#1{\mathord{\mathop{\kern 0pt #1}\limits^\circ}}

\title[]{Equivariant isometric immersions of surfaces in hyperbolic space}

\author{Jean-Marc Schlenker}
\address{Jean-Marc Schlenker:
University of Luxembourg, FSTM, Department of Mathematics,
Maison du nombre, 6 avenue de la Fonte,
L-4364 Esch-sur-Alzette, Luxembourg}
\email{jean-marc.schlenker@uni.lu}
\thanks{Partially supported by FNR project O25/19559635/SIDE.}

\begin{document}

\begin{abstract}
  Let $S$ be a closed, oriented surface of genus at least $2$, let $h$ be a smooth Riemannian metric on $S$ with curvature $K\in (-1,0]$, and let $c\in \cT_S$ be a conformal structure on $S$. There exists a unique equivariant isometric immersion of $(S,h)$ in $\HH^3$ such that the pull-back by the hyperbolic Gauss map of the conformal structure at infinity is $c$. Dually, if $h^*$ is a smooth metric on $S$ with curvature $K^*\in (-\infty, 0)$ and if $c\in \cT_S$, there exists a unique equivariant immersion of $S$ into $\HH^3$ with third fundamental form $h^*$ and such that the pull-back of the conformal class at infinity by the Gauss map is $c$.

  Equivalently, given $h$ and $c\in \cT_S$, there is a unique pair $(E,u)$ where $E$ is a hyperbolic end with conformal structure at infinity $c$ and $u$ is an isometric embedding of $(S,h)$ in $E$. Given $h^*$ and $c$, there exists a unique pair $(E,u^*)$, where $E$ is a hyperbolic end with conformal structure at infinity $c$, and $u^*:S\to E$ is an embedding inducing the third fundamental form $h^*$.

  Those statements can be considered as a smooth counterpart of known or conjectural statements on the grafting map and on circle patterns.
\end{abstract}

\maketitle

\tableofcontents

\section{Results and motivations}

\subsection{The hyperbolic Weyl problem and its dual}

The Weyl problem has been a driving force in the development of differential geometry in the twentieth century. Its hyperbolic version, solved by Alexandrov \cite{alex} and Pogorelov \cite{Po},  asks whether, given a Riemannian metric $h$ on $S^2$ with curvature $K>-1$, it admits a unique isometric immersion in $\HH^3$. The curvature condition then implies that the image surface is locally convex, which in turn implies that it bounds a convex subset in $\HH^3$. 

A dual question concerns ``bending'' invariants of the boundary, such as the third fundamental form. In the polyhedral setting, it corresponds to a result of Hodgson--Rivin \cite{HR} describing the possible ``dual metrics'' (which are closely related to the dihedral angles) of compact polyhedra in $\HH^3$, while a smooth statement is in \cite{these}.

Our goal here is to extend those results to surfaces of higher genus, specifically of genus at least $2$. For such a higher genus surface $S$,  one can consider {\em equivariant} immersions of $S$ in $\HH^3$, that is, pairs $(u, \rho)$, where $u:\tilde S\to \HH^3$ is an immersion of the universal cover of $S$, $\rho:\pi_1S\to \Isom(\HH^3)$ is a morphism, and
$$ \forall x\in \tilde S, \forall \gamma\in \pi_1S, u(\gamma.x)=\rho(\gamma)(u(x))~. $$

Those equivariant immersions, however, are too ``flexible'' to hope for a uniqueness result when prescribing the induced metric (or, for the dual statement, the third fundamental form). In fact, a dimension count shows that the space of equivariant immersions of $S$ into $\HH^3$ with prescribed induced metric should have dimension $6g-6$, where $g$ is the genus of $S$. So additional data needs to be prescribed --- clearly, this doesn't happen for spheres, in the classical Weyl problem.

Here we propose that this additional data is the conformal structure at infinity associated to a locally convex equivariant immersion of $S$, as described below. We then prove an existence and uniqueness result for equivariant isometric immersions of $(S,h)$ into $\HH^3$, when $h$ has curvature $K\in (-1,0]$. We also prove a dual existence and uniqueness statement, for equivariant immersions with prescribed third fundamental form $h^*$, when $h^*$ has curvature $K^*< 0$.

There is another, equivalent way to state the hyperbolic Weyl problem and its dual, in terms of the induced metric or third fundamental form on the boundary of convex domains in $\HH^3$. This leads to a different type of generalizations, for instance on the induced metric or third fundamental form (or bending measure, etc) on the boundary of convex subsets in hyperbolic manifolds, see \cite{weylsurvey}.

\subsection{Equivariant isometric immersions}

We first define the notions appearing in the main results. Here and below we denote by $S$ a closed, oriented surface, which will later need to have genus at least $2$. 

\begin{definition}
  \begin{itemize}
  \item An {\em equivariant immersion} of $S$ in $\HH^3$ is a pair $(u,\rho)$ where $u:\tilde S\to \HH^3$ is an immersion, $\rho:\pi_1S\to \Isom(\HH^3)$, and
    $$ \forall x\in \tilde S, \forall \gamma\in \pi_1S, u(\gamma.x)= \rho(\gamma)(u(x))~. $$
  \item $(u,\rho)$ is {\em locally convex} if the second fundamental form of $u(\tilde S)$ is positive definite.
  \item The hyperbolic Gauss map of $u$ is the map $G_u:\tilde S\to \partial_\infty\HH^3$ sending $x\in \tilde S$ to the endpoint in $\partial_\infty \HH^3$ of the geodesic ray starting from $u(x)$ in the direction of the oriented unit normal to $u(\tilde S)$ at $u(x)$.
  \end{itemize}
\end{definition}

Here the second fundamental form of $u(\tilde S)$ is defined on two vectors $v,w\in T_x\tilde S$ by
$$ \II(v,w) = g(D_vN,w)~, $$
where $N$ is the oriented unit normal to $u(\tilde S)$ at $u(x)$, $g$ is the Riemannian metric of $\HH^3$, and $D$ is its Levi-Civita connection.

\begin{theorem} \label{tm:main}
  Let $S$ be a closed, oriented surface of genus at least $2$, let $h$ be a smooth Riemannian metric on $S$ with curvature $K\in (-1,0]$, and let $c\in \cT_S$ be a conformal structure on $S$. There exists a unique locally convex equivariant isometric immersion of $(S,h)$ in $\HH^3$ such that the pull-back by the hyperbolic Gauss map of the conformal structure at infinity is isotopic to $c$.
\end{theorem}

This is a positive answer, for the special case where $K\in (-1,0]$, of \cite[Question $W^{imm}_{\HH^3}$]{weylsurvey}, which is stated there for $K>-1$. Uniqueness is up to composition by a global isometry of $\HH^3$.

The dual statement also holds, for surfaces of negative curvature. For surfaces of genus $0$ it reduces to the main result of \cite{these}, but again we need to assume, for surfaces of higher genus, that the curvature is negative.

\begin{theorem} \label{tm:main*}
  Let $S$ be a closed, oriented surface of genus at least $2$, let $h^*$ be a smooth Riemannian metric on $S$ with curvature $K^*<0$, and let $c\in \cT_S$ be a conformal structure on $S$. There exists a unique locally convex equivariant immersion of $S$ in $\HH^3$ with third fundamental form $h^*$ such that the pull-back by the hyperbolic Gauss map of the conformal structure at infinity is isotopic to $c$.
\end{theorem}

Recall that the third fundamental form of an immersion is defined, with the notations above, as
$$ \III(v,w) = g(D_vN,D_wN)~. $$
This is again the special case of \cite[Question $W^{imm}_{\HH^3}$]{weylsurvey} where we suppose that $K^*$, the curvature of $h^*$, is negative. (The question as stated in \cite{weylsurvey} also has a condition on the lengths of closed, contractible geodesics, which is vacuous when $K^*<0$.) Here too, uniqueness is up to composition by a global isometry of $\HH^3$.

\subsection{Hyperbolic ends} \label{ssc:ends}

We now recall the definition of a hyperbolic end used here.

\begin{definition}\label{df:end}
  Let $S$ be a closed oriented surface of genus at least $2$. A {\em hyperbolic end} is a hyperbolic $3$-manifold $E$, diffeomorphic to $S\times\R_{>0}$, such that:
  \begin{itemize}
  \item the metric is non-complete, and complete on the side corresponding to $\infty$,
  \item the metric completion $\bar E$ is obtained by adding a pleated surface $\partial_0E$ corresponding to $S\times\{0\}$,
  \item $E$ is concave in the neighbourhood of $\partial_0E$, that is, the universal cover $\tilde E$ is locally isometric to the complement in $\HH^3$ of a closed convex set with non-empty interior.
  \end{itemize}
  We call $\partial_0E$ the {\em pleated boundary} of $E$ and $\partial_\infty E$ its {\em boundary at infinity}, and we denote by $\cE_S$ the space of hyperbolic ends diffeomorphic to $S\times\R_{>0}$, considered up to isotopy.
\end{definition}

We fix the following additional conventions. The metric completion of $E$ is $E\cup\partial_0 E$, where $\partial_0 E$ is a locally concave pleated surface (where the concave side is $E$). Its induced path metric is a hyperbolic metric $m(E)\in\cT_S$ and its bending is described by a measured lamination $l(E)\in\cML_S$, its {\em bending lamination}. The boundary at infinity $\partial_\infty E$ carries a $\bCP^1$-structure $\sigma(E)\in\cP_S$ whose underlying conformal structure is $c(E)\in\cT_S$. Both $\partial_0 E$ and $\partial_\infty E$ are identified with $S$ up to isotopy. For $(E,u)\in\cE$ the hyperbolic Gauss map $u(S)\to\partial_\infty E$ is isotopic to this identification, so $c(E)$ is the conformal structure appearing in Theorem \ref{tm:main}.

A connected component of the complement of the convex core of a convex co-compact hyperbolic manifold is a hyperbolic end, but the converse does not hold, since a hyperbolic end need not embed in a complete hyperbolic manifold. The definition above excludes degenerate ends --- the terminology ``convex co-compact hyperbolic end'' would be more precise, see \cite[Def. 1.3]{cp}.

Given an oriented surface $S\subset E$, its hyperbolic Gauss map is the map sending $x\in S$ to the end point of the geodesic ray starting from $x$ in the direction of the oriented normal of $S$ at $x$. This endpoint could be in $\partial_\infty E$, in $\partial_0 E$, or it might be undefined. However we are only interested in the hyperbolic Gauss map for surfaces for which it takes values in $\partial_\infty E$. 

\begin{definition}
  Let $E$ be a hyperbolic end, and let $u:S\to E$. We will say that $u$ is a {\em good embedding} in $E$ if its image is locally convex and the hyperbolic Gauss map is a diffeomorphism from $u(S)$ to $\partial_\infty E$. We will say that a surface in $E$ is  a {\em good surface} if it is the image of a good embedding.
\end{definition}

If $u$ is a good embedding in $E$, then $u(S)$ is locally {\em strongly convex}, that is, its principal curvatures are strictly positive. It then follows from the Gauss formula that its induced metric has curvature $K>-1$.

For instance, the boundary of an $r$-neighborhood of $\partial_0E$ is a good surface. However, if $x\in E$, the boundary of the ball of center $x$ and radius $r$ is not (for $r$ small enough) a good surface.

\begin{lemma} \label{lm:equiv}
  Let $(u,\rho)$ be a locally convex equivariant immersion of $S$ in $\HH^3$. Then there exists a unique hyperbolic end $E$ and a unique good embedding $v:S\to E$ such that $v$ lifts to $u$. Conversely, a good embedding $v:S\to E$ in a hyperbolic end lifts to an equivariant locally convex immersion of $S$ in $\HH^3$. 
\end{lemma}

The proof is in Section \ref{ssc:conv-ends}.

We can therefore restate Theorems \ref{tm:main} and \ref{tm:main*} in an equivalent form as describing locally convex embeddings in hyperbolic ends.

\begin{corollary} \label{cr:main}
  Let $h$ be a smooth metric on $S$ of curvature $K\in (-1,0]$, and let $c\in \cT_S$ be a conformal structure on $S$. There is a unique pair $(E,u)$ where $E$ is a hyperbolic end with conformal structure at infinity isotopic to $c$ and $u$ is a good isometric embedding of $(S,h)$ in $E$.
\end{corollary}

\begin{corollary} \label{cr:main*}
    Let $h^*$ be a smooth metric on $S$ of curvature $K^*<0$, and let $c\in \cT_S$ be a conformal structure on $S$. There is a unique pair $(E,u)$ where $E$ is a hyperbolic end with conformal structure at infinity isotopic to $c$ and $u$ is a good embedding of $S$ in $E$ with third fundamental form $h^*$.
\end{corollary}

\subsection{Smooth grafting} \label{ssc:sgr}

The smooth grafting map $sgr':\R_{>0}\times \cT_S\times \cT_S\to \cT_S$ is defined in \cite{cyclic} (or \cite[\S 2.10]{cyclic2}). It can be described as follows. Let $s>0$, and let $h,h^*\in\cT_S$. There is then a unique hyperbolic end $E$ containing an embedded surface with induced metric $\cosh^2(s/2)h$ and third fundamental form $\sinh^2(s/2)h^*$. Then $sgr'(s,h,h^*)$ is the conformal class at infinity of $E$. The smooth grafting map can be seen as an ``imaginary counterpart'' to the landslide map of \cite{cyclic}, in the sense that they fit together to make up a ``complex landslide'' map.

It is proved in \cite[Theorem 1.7]{cyclic2} that for all $s>0$, the maps $sgr'(s,h,\cdot):\cT_S\to \cT_S$ and $sgr'(s,\cdot,h^*):\cT_S\to \cT_S$ are proper and surjective. It was then asked whether they are injective. As a direct consequence of Corollary \ref{cr:main} and Corollary \ref{cr:main*} we obtain:

\begin{corollary} \label{cr:sgr}
  For all $s>0$ and all $h,h^*\in \cT_S$, the maps $sgr'(s,h,\cdot):\cT_S\to \cT_S$ and $sgr'(s,\cdot,h^*):\cT_S\to \cT_S$ are homeomorphisms.
\end{corollary}

\subsection{Dual surfaces in de Sitter ends} \label{ssc:dual}

Theorems \ref{tm:main} and \ref{tm:main*} can also be stated, dually, in terms of surfaces embedded in a globally hyperbolic compact maximal (GHMC) de Sitter spacetime.

To each hyperbolic end, one can associate a GHMC de Sitter spacetime. One way to see this is by noting that hyperbolic ends homeomorphic to $S\times \R_{>0}$ are in one-to-one correspondence with complex projective structures on $S$ (see e.g. \cite{kulkarni-pinkall}), while complex projective structures on $S$ are also in one-to-one correspondence with GHMC de Sitter spacetimes \cite{mess,scannell}.

The correspondence is easier to visualize when the developing map of the hyperbolic end $E$ is injective (or equivalently the developing map of its complex projective structure at infinity is injective). In this case $E=\tilde E/\rho(\pi_1S)$, where $\tilde E\subset \HH^3$. We can define the dual domain $\tilde E^*\subset dS^3$ as the space of points in $dS^3$ which are dual to oriented planes in $\HH^3$ which are contained in $\tilde E$ and bounding a half-space containing $\widetilde{\partial_0E}$. It is invariant under the action of $\pi_1S$ and the quotient, $E^*=\tilde E^*/\rho(\pi_1S)$, is a GHMC de Sitter spacetime. The same construction can in fact be done (in a somewhat less visual manner) when the developing map of $E$ is not injective, and points of $E^*$ are still dual to the oriented planes immersed in $E$, in the sense of the hyperbolic--de Sitter duality (see e.g. \cite{shu,fillastre-seppi}). The complex projective structure at infinity of $E$ is the same as that of $E^*$, and the same therefore also applies to the conformal structure at infinity.

If $u:S\to E$ is a good embedding, one can define a dual embedding $u^*:S\to E^*$. The metric induced on $S$ by $u$ is the third fundamental form induced by $u^*$, and conversely. Theorems \ref{tm:main} and \ref{tm:main*} can therefore be stated in the following, equivalent forms:
\begin{itemize}
\item Let $h$ be a smooth metric on $S$ of curvature $K\in (-1,0]$, and let $c\in \cT_S$ be a conformal structure on $S$. There is a unique pair $(E^*,u^*)$ where $E^*$ is a GHMC de Sitter spacetime with conformal structure at infinity isotopic to $c$ and $u^*$ is a spacelike embedding of $S$ in $E^*$ with third fundamental form $h$.
\item Let $h^*$ be a smooth metric on $S$ of curvature $K^*< 0$, and let $c\in \cT_S$ be a conformal structure on $S$. There is a unique pair $(E^*,u)$ where $E^*$ is a GHMC de Sitter spacetime with conformal structure at infinity isotopic to $c$ and $u$ is a spacelike embedding of $S$ in $E^*$ with induced metric $h^*$.
\end{itemize}

\subsection{Motivations and limit cases}

There are a number of relatively recent results which can be seen as either special cases or limit cases of Theorems \ref{tm:main} and \ref{tm:main*}, or which have close connections to them (in addition to the applications to smooth grafting in Section \ref{ssc:sgr}). Moreover, an interesting question on circle packings, or circle patterns, can also be considered as a limit case of Theorem \ref{tm:main*}. 

For metrics of constant curvature, the existence part of Theorem \ref{tm:main} was proved in \cite[Thm 1.9]{cyclic2}, where the question was also raised of the uniqueness. So the uniqueness part of Theorem \ref{tm:main}, for metrics of constant curvature, can be stated in terms of the landslides studied in \cite{cyclic2} --- we leave the details to the interested reader.

In the limit case where $K=-1$, Theorem \ref{tm:main} corresponds to a result of Dumas and Wolf \cite{dumas-wolf}. They proved that the grafting map for a fixed hyperbolic metric $m$ on $S$, considered as a map $gr_m:\cML_S\to \cT_S$, is a homeomorphism. 

Still in the limit case where $K=-1$, or as $K^*\to -\infty$, Theorem \ref{tm:main*} corresponds to a result of Scannell and Wolf \cite{scannell-wolf}. In this limit case the third fundamental form of smooth surfaces corresponds --- in a limit sense --- to the measured bending lamination of pleated surfaces. Scannell and Wolf proved that, for a fixed measured lamination, the grafting map $gr_l:\cT_S\to \cT_S$, is a homeomorphism.

Another limit case of Theorem \ref{tm:main*} concerns ideal polyhedral surfaces. Recall that if $\Sigma$ is a convex polyhedral surface in a hyperbolic end $E$, its third fundamental form still makes sense as the transverse measure on its edges, giving to each edge a weight equal to its exterior dihedral angle.
If  all the vertices of $\Sigma$ are ideal, the link of an ideal vertex $v$ is a Euclidean polygon, so the face angles at $v$ add up to $2\pi$. The dual metric
of an ideal polyhedral surface therefore corresponds to the graph $\Gamma^*$ dual to the 1-skeleton of the polyhedral surface, with edge $e$ of length $\theta_e$. So Theorem \ref{tm:main*} (in a slightly more general form, requesting only that $K^*<1$) appears as a conjecture of Kojima, Mizushima and Tan \cite{KMT,KMT2,KMT3} on circle packings, in the more general form appearing in \cite{delaunay}, where circle packings are replaced by Delaunay circle patterns.

\begin{question} \label{q:KMT}
  Let $\Gamma$ be a graph embedded in $S$, which is the 1-skeleton of a cell decomposition. Let $\theta:E(\Gamma)\to (0,\pi)$ be such that
  \begin{enumerate}
  \item the sum of values of $\theta$ at each vertex is equal to $2\pi$,
  \item the sum of values of $\theta$ on any homotopically trivial non-backtracking closed edge path $[e_1,\cdots, e_n]$ in the dual graph $\Gamma^*$ which does not bound a face is larger than $2\pi$.
  \end{enumerate}
  Finally let $c\in \cT_S$. Then there exists a unique pair $(E,\Sigma)$, where $E$ is a hyperbolic end and $\Sigma$ is an ideal polyhedral surface in $E$ with combinatorics $\Gamma$ and exterior dihedral angles $\theta$, and such that the conformal structure at infinity of $E$ is $c$.
\end{question}

Condition (1) in Question \ref{q:KMT} states that each face of $\Gamma^*$ has boundary length exactly $2\pi$, which is precisely what forces the vertices of $\Sigma$ to be ideal. Condition (2) is the condition, which would be necessary in Corollary \ref{cr:main*} if we didn't restrict the curvature to be negative, that all closed contractible geodesics have length larger than $2\pi$ --- the boundaries of the faces of $\Gamma^*$ being excluded since they are exactly the closed geodesics realizing the limit value $2\pi$.

So Question \ref{q:KMT} is a limit case of the more general statement, for metrics of curvature $K^*<1$ all of whose closed contractible geodesics have length larger than $2\pi$, which is \cite[Question $W^{imm}_{\HH^3}$]{weylsurvey}, and which we do not prove here.

In a different direction, there is another way to obtain an existence and uniqueness statement: rather than by prescribing the pull-back of the conformal structure at infinity by the Gauss map, one can request that the hyperbolic end (or equivalently the holonomy representation of the equivariant immersion) is Fuchsian, that is, contained in $\PSL(2,\R)\subset \PSL(2,\C)$. Given a metric $h$ on $S$ of curvature $K>-1$, there is a unique such Fuchsian equivariant immersion in $\HH^3$ with induced metric $h$ \cite[3.2.4 B', B'']{PDR},\cite{fillastre2,hmcb}. Dually, given a smooth metric $h^*$ on $S$ of curvature $K^*<1$, with no closed, contractible geodesic of length $L\leq 2\pi$, there exists a unique Fuchsian equivariant immersion of $S$ in $\HH^3$ inducing the third fundamental form $h^*$, see \cite{iie} (via duality).

\subsection{Sketch of the proofs} \label{ssc:sketch}

The proof uses a variation on a deformation argument. We introduce some notations first. 

\begin{definition}
  \begin{itemize}
  \item $\cM_S$ is the space of smooth metrics on $S$ of curvature $K\in (-1,0]$, and  $\bar \cM_S$ is the space of smooth metrics on $S$ of curvature $K>-1$.
  \item $\cM_S^*$ is the space of smooth metrics on $S$ of curvature $K\in (-\infty,0)$, and $\bar \cM_S^*$ is the space of smooth metrics on $S$ of curvature $K<1$.
  \item $\cE$ is the space of pairs $(E,u)$, where $E$ is a hyperbolic end and $u:S\to E$ is a good embedding of $S$ in $E$, isotopic to the canonical identification of $S$ with $S\times \{1\}\subset E$.
  \item For all $h\in \bar \cM_S$, $\cG_h$ is the space of solutions of the Gauss-Codazzi equations on $(S,h)$, that is, the space of $h$-self-adjoint positive definite bundle morphisms $B:TS\to TS$ such that $K_h=-1+\det(B)$ and $d^DB=0$, where $D$ is the Levi-Civita connection of $h$.
  \item For all $h^*\in \bar \cM_S^*$, $\cG^*_{h^*}$ is the space of solutions of the dual Gauss-Codazzi equations on $(S,h^*)$, that is, the space of $h^*$-self-adjoint positive definite bundle morphisms $B^*:TS\to TS$ such that $K_{h^*}=1-\det(B^*)$ and $d^{D^*} B^*=0$, where $D^*$ is the Levi-Civita connection of $h^*$.
  \item $\cm:\cE\to \bar \cM_S$ is the map sending a pair $(E,u)$ to the metric induced on $S$ by $u$.
  \item $\ct:\cE\to \bar \cM_S^*$ is the map sending a pair $(E,u)$ to the pull-back by $u$ of the third fundamental form of $u(S)$.
  \item $\cc:\cE\to \cT_S$ is the map sending a pair $(E,u)$ to the conformal class at infinity of $E$.
  \item For all $h\in \cM_S$, $\cc_h:\cG_h\to \cT_S$ is the map sending a solution $B$ of the Gauss-Codazzi equations on $(S,h)$ to the conformal class at infinity of the hyperbolic end containing the surface with induced metric $h$ and shape operator $B$.
  \item For all $h^*\in \cM^*_S$, $\cc^*_{h^*}:\cG^*_{h^*}\to \cT_S$ is the map sending a solution $B^*$ of the dual Gauss-Codazzi equations on $(S,h^*)$ to the conformal class at infinity of the hyperbolic end containing the surface with third fundamental form $h^*$ and shape operator $B=(B^*)^{-1}$.
  \end{itemize}
\end{definition}

\subsubsection{Sketch of the proof of Theorem \ref{tm:main}}

We first give an outline of the proof of Theorem \ref{tm:main}. As mentioned already, it relies on a deformation argument: we need to prove that for all $h\in \cM_S$, the map $\cc_h:\cG_h\to \cT_S$ is a homeomorphism. The first point is that this map is a local homeomorphism. It will follow from the fact that both the domain and the target have dimension $6g-6$, and the following lemma.

\begin{lemma}\label{lm:rigid}
Let $\Sigma$ be a good surface in a hyperbolic end, with induced metric of curvature $K\in (-1,0]$. Then every first-order deformation of the pair $(E,\Sigma)$ preserving the induced metric on $\Sigma$ and the isotopy class of the conformal class at infinity of $E$ is trivial.
\end{lemma}

The second step is to prove that $\cc_h$ is in fact proper, which implies that it is a covering map.

\begin{lemma} \label{lm:proper}
  Let $(h_n)_{n\in \N}$ be a sequence of smooth metrics on $S$ with curvature $K>-1$, converging smoothly to a smooth metric $h$ with $K>-1$. Let $(E_n,u_n)_{n\in \N}$ be a sequence of pairs where $E_n$ is a hyperbolic end with conformal structure at infinity $c_n$ and $u_n:(S, h_n)\to E_n$ is a good isometric embedding. Assume that $c_n\to c\in \cT_S$. Then, after extracting a subsequence, $(E_n, u_n)_{n\in \N}$ converges smoothly to a limit $(E,u)$, where $E$ is a hyperbolic end with conformal structure at infinity $c$ and $u:(S,h)\to E$ is a good isometric embedding of $(S,h)$ in $E$.
\end{lemma}

Theorem \ref{tm:main} then follows from the fact that $\cT_S$ is simply connected and, when $h$ has constant curvature, that $\cG_h$ is connected, which implies that any covering map from $\cG_h$ to $\cT_S$ has degree $1$ and is therefore a homeomorphism. 

\subsubsection{Sketch of the proof of Theorem \ref{tm:main*}}

The same argument can be used, using the map $\cc^*_{h^*}:\cG^*_{h^*}\to \cT_S$. The infinitesimal rigidity lemma is now the following statement.

\begin{lemma}\label{lm:rigid*}
Let $\Sigma$ be a good surface in a hyperbolic end, with induced metric of curvature $K\in (-1,0]$. Then every first-order deformation of the pair $(E,\Sigma)$ preserving the third fundamental form on $\Sigma$ and the isotopy class of the conformal class at infinity of $E$ is trivial.
\end{lemma}

The properness of $\cc^*_{h^*}$ corresponds to the following lemma.

\begin{lemma} \label{lm:proper*}
  Let $(h_n^*)_{n\in \N}$ be a sequence of smooth metrics on $S$ with curvature $K^*<0$, converging smoothly to a smooth metric $h^*$ with $K^*<0$. Let $(E_n,u^*_n)_{n\in \N}$ be a sequence of pairs where $E_n$ is a hyperbolic end with conformal structure at infinity $c_n$ and $u^*_n:S\to E_n$ is a good embedding with third fundamental form $h_n^*$. Assume that $c_n\to c\in \cT_S$. Then, after extracting a subsequence, $(E_n, u^*_n)_{n\in \N}$ converges smoothly to a limit $(E,u^*)$, where $E$ is a hyperbolic end with conformal structure at infinity $c$ and $u^*:S\to E$ is a good embedding of $S$ in $E$ inducing the third fundamental form $h^*$.
\end{lemma}

The same topological argument as above then proves the result.

\subsection{Note on use of AI}

This manuscript was written thanks to a substantial help from AI tools, specifically Claude Opus 5, Fable 5 and Fable 5.1. Specifically, Fable 5, after being provided with pages of notes of the author, made a key contribution to the proof of Lemma \ref{lm:rigid}, extending slightly, but in a very significant way, on the computations and argument provided in the notes.

A first draft of the proofs of Lemma \ref{lm:proper} and Lemma \ref{lm:proper*} was also written by Fable 5 and 5.1, but only after many rounds of increasingly precise instructions from the author. It appears likely that writing the proof of those lemmas directly without the help of AI might have been faster. However even for this part, AI was helpful, for instance in pointing out a gap in an earlier version of an argument. Other models, such as Opus 5 and Opus 5.5, were used for lighter tasks, such as proofreading.

The final version of the manuscript was heavily human-edited, since the presentations in the AI-produced text didn't seem suitable to the author. 

\subsection{Acknowledgement}

I would like to thank Francesco Bonsante, Qiyu Chen and Wayne Lam for interesting conversations or email exchanges related to this project.

\section{Background}

\subsection{Conventions}\label{ssc:conventions}

Throughout, $S$ is closed, oriented, of genus $g\geq2$. If $I$ is a Riemannian metric on $S$, we write $D$ for its Levi-Civita connection, $K$ for its curvature, $da$ for its area form, and $J$ for the rotation by $+\pi/2$; thus $R(x,y)z=D_xD_yz-D_yD_xz-D_{[x,y]}z$ satisfies $R(x,y)z=-K\,da(x,y)\,Jz$.

For a bundle endomorphism $A:TS\to TS$, we write
$$ d^DA(x,y)=D_x(Ay)-D_y(Ax)-A[x,y]~. $$
$A$ is \emph{Codazzi} if $d^DA=0$. On an oriented $2$--plane, $\adj(A)=(\tr A)E-A$ denotes the adjugate, so that $\adj(A)A=(\det A)E$ and, for any $A$,
\begin{equation}\label{eq:adj}
\adj(A)=-JAJ,\qquad \tr(\adj(A)B)=\tr A\tr B-\tr(AB).
\end{equation}

We will repeatedly use that $\tr(JA)=0$ if and only if $A$ is self-adjoint, and that $(JA)^2=-(\det A)E$ for $A$ self-adjoint.

\subsection{Identities on $2\times 2$ matrices}

We collect the pointwise linear algebra used below. For an endomorphism $A$ of a $2$-dimensional vector space, $\adj(A)A=A\adj(A)=(\det A)E$ by the Cayley--Hamilton theorem, and $A\mapsto\adj(A)$ is linear.

\begin{lemma}\label{lm:2d}
  Let $(V,\langle\cdot,\cdot\rangle)$ be an oriented Euclidean plane and let
  $J$ be the rotation by $+\pi/2$. For all self-adjoint endomorphisms $P,Q$ of
  $V$:
  \begin{enumerate}
  \item $\adj(P)=-JPJ$, and consequently $J\adj(P)=PJ$;
  \item $\tr(\adj(P)Q)=\tr P\,\tr Q-\tr(PQ)=\tr(\adj(Q)P)$;
  \item $\tr(JP)=0$;
  \item $\tr(JPQ)=-\tr(JQP)$;
  \item $\tr(JPJP)=-2\det P$ and $\tr(\adj(P)JP)=0$.
  \end{enumerate}
\end{lemma}

\begin{proof}
  (1) is a direct verification in a positively oriented orthonormal basis, in
  which $P$ is represented by a symmetric matrix; then
  $J\adj(P)=-J^2PJ=PJ$. (2) is immediate from the definition of $\adj$, and
  manifestly symmetric in $P,Q$. For (3), taking adjoints,
  $\tr(JP)=\tr((JP)^t)=\tr(P(-J))=-\tr(JP)$. For (4), similarly
  $\tr(JPQ)=\tr((JPQ)^t)=\tr(QP(-J))=-\tr(JQP)$. For (5), by (1),
  $JPJP=(JPJ)P=-\adj(P)P=-(\det P)E$, whose trace is $-2\det P$, and
  $\adj(P)JP=(-JPJ)JP=JP^2$, whose trace vanishes by (3) since $P^2$ is
  self-adjoint.
\end{proof}

\subsection{Divergence and Codazzi operators}

In this section we recall basic facts concerning Codazzi operators on surfaces, for later use. To keep notations simple, we still consider the closed surface $S$, equipped with a Riemannian metric which we call $I$, since it will later be the induced metric (or in some cases the third fundamental form) of an embedding.

\begin{lemma}\label{lm:div}
  Let $A$ be a smooth field of $I$-self-adjoint endomorphisms of $TS$ and let
  $u$ be a vector field on $S$. Set $\dv A=\sum_i(D_{e_i}A)e_i$ for a local
  orthonormal frame $(e_i)_{i=1,2}$, and write $Du$ for the endomorphism $x\mapsto
  D_xu$. Then:
  \begin{enumerate}
  \item $\dv(Au)=I(u,\dv A)+\tr(A\,Du)$.
  \item If moreover $d^DA=0$, then $\dv A=D(\tr A)$ and $\dv(\adj(A))=0$. In
    particular, since $S$ is closed,
    $$ \int_S \tr(\adj(A)\,Du)\,da=0 $$
    for every vector field $u$.
  \item If $B:TS\to TS$ is a bundle morphism which is self-adjoint for $I$ and Codazzi, then $\dv(J(E+B)u)=\tr(J(E+B)\,Du)$.
  \end{enumerate}
\end{lemma}

\begin{proof}
  All computations are done at a point where $De_i=0$. Note first that $D_XA$
  is self-adjoint for every $X$, as follows by differentiating
  $I(Ax,y)=I(x,Ay)$ along $X$.

  (1) We have
  $$ \dv(Au)=\sum_i I(D_{e_i}(Au),e_i)
  =\sum_i I((D_{e_i}A)u,e_i)+\sum_i I(A\,D_{e_i}u,e_i)~. $$
  The second sum is $\tr(A\,Du)$ by definition, and the first equals
  $\sum_i I(u,(D_{e_i}A)e_i)=I(u,\dv A)$ by self-adjointness of $D_{e_i}A$.

  (2) For each $j$, using successively the self-adjointness of $D_{e_i}A$ and
  the Codazzi equation $(D_{e_i}A)e_j=(D_{e_j}A)e_i$, if $u$ is a vector field on $S$,
  $$ I(\dv A,u)=\sum_i I((D_{e_i}A)e_i,u)=\sum_i I(e_i,(D_{e_i}A)u)
  =\sum_i I(e_i,(D_{u}A)e_i)=\tr(D_{u}A)=d(\tr A)(u)~, $$
  so $\dv A=D(\tr A)$. Since $\dv(fE)=Df$ for a function $f$, it follows that
  $$ \dv(\adj(A))=\dv((\tr A)E)-\dv A=0~. $$
  The adjugate $\adj(A)$ is self-adjoint, so by (1) applied to $\adj(A)$,   $\dv(\adj(A)u)=\tr(\adj(A)\,Du)$, and integrating over the closed surface $S$ gives the last assertion.

  (3) Since $J$ and $E$ are parallel,
  $$ \dv(J(E+B)u)=\sum_i I(J(D_{e_i}B)u,e_i)+\tr(J(E+B)\,Du)~. $$
  By the Codazzi equation, $(D_{e_i}B)u=(D_uB)e_i$, so the first sum equals $\sum_i I(J(D_uB)e_i,e_i)=\tr(J\,D_uB)$, which vanishes by Lemma~\ref{lm:2d}(3) since $D_uB$ is self-adjoint.
\end{proof}

\subsection{Locally convex equivariant isometric immersions and hyperbolic ends} \label{ssc:conv-ends}

A {\em hyperbolic end} is a non-complete 3-dimensional hyperbolic manifold $E$ associated to a complex projective structure $\sigma$ on a closed surface $S$ of genus at least $2$. An example is provided by each connected component of the complement of the convex core in a quasifuchsian (or more generally a convex co-compact) hyperbolic manifold: $E$ is homeomorphic to $S\times \R_{>0}$, complete on the side of $+\infty$, but non-complete on the side corresponding to $0$, where it is bounded by a concave pleated surface. In this case, and when the developing map of the complex projective structure $\sigma$ on $S$ is injective, the hyperbolic end $E$ can be defined as
$$ E = (\HH^3\setminus CH(\CP^1\setminus \mathrm{dev}_\sigma(\tilde S)))/\rho_\sigma(\pi_1S)~, $$
where $CH$ denotes the convex hull, and $\mathrm{dev}_\sigma$ and $\rho_\sigma$ are the developing map and holonomy representation of $\sigma$.

When the developing map of $\sigma$ is not injective, the construction of $E$ is less visual, and involves the gluing of half-spaces associated to maximal disks in $(\tilde S,\sigma)$, see \cite{kulkarni-pinkall,kamishima-tan}.

\begin{proof}[Proof of Lemma \ref{lm:equiv}]
  Let $E$ be a hyperbolic end, and let $v:S\to E$ be a good embedding. Then, by definition, $v$ lifts to an equivariant locally convex immersion $u:\tilde S\to \HH^3$, with holonomy map $\rho$ equal to the holonomy representation of $E$.

  Conversely, let $u:\tilde S\to \HH^3$ be a locally convex equivariant immersion of $S$ into $\HH^3$, with holonomy map $\rho:\pi_1S\to \mathrm{Isom}(\HH^3)$. Let $U: \tilde S\times \R_{\geq 0}\to \HH^3$ be the map sending $(x,t)$ to $\gamma_x(t)$, where $\gamma_x:\R_{\geq 0}\to \HH^3$ is the unit-velocity parameterized geodesic ray starting from $u(x)$ in the direction of the oriented normal to $du(T_xS)$. Since $u$ is locally convex, $U$ is an immersion, and the pull-back by $U$ of the hyperbolic metric defines a hyperbolic metric $g$ on $\tilde S\times \R_{\geq 0}$, which is non-complete on $\tilde{S}\times \{ 0\}$.

  By construction, the action $\rho$ of $\pi_1S$ on $\tilde S$ extends to an isometric action (that we still denote by $\rho$) of $\pi_1S$ on $(\tilde S\times \R_{\geq 0}, g)$. The quotient is a hyperbolic manifold $E_0$, naturally identified with $S\times \R_{\geq 0}$, which is non-complete on $S\times \{ 0\}$, but equips $S\times \{\infty\}$ with a $\CP^1$-structure $\sigma$.

  Still by construction, the boundary $S\times \{ 0\}$ is concave. As a consequence, $E_0$ embeds isometrically in a unique hyperbolic end $E$, which has the $\CP^1$-structure $\sigma$ at infinity, see \cite{kulkarni-pinkall}, or more specifically \cite[\S 3.8]{grothu}.

  Still by construction, $u$ determines a good embedding $v:S\to E$, with image $S\times \{ 0\}\subset E_0\subset E$. 
\end{proof}

\begin{proof}[Proof of Corollary \ref{cr:main} and of Corollary \ref{cr:main*}]
  Both corollaries follow from the corresponding theorem (Theorem \ref{tm:main} and Theorem \ref{tm:main*}, respectively) because, if $u:\tilde S\to \HH^3$ is a locally convex equivariant immersion of $S$ in $\HH^3$, and if $v:S\to E$ is the corresponding good embedding of $S$ in a hyperbolic end, then the pull-back on $S$ of the conformal class at infinity of $\HH^3$ by the Gauss map of $v$ is isotopic to the conformal class at infinity of $E$. 
\end{proof}

\subsection{Surfaces in $\HH^3$}

We will use the following notations on oriented surfaces in $\HH^3$.
\begin{itemize}
\item $I$ denotes the induced metric on $S$, that is, the restriction to $TS$ of the ambient metric.
\item $N$ is the oriented unit normal vector field to $S$.
\item $D$ is the Levi-Civita connection of $\HH^3$.
\item $B:TS\to TS$ is the shape operator of $S$, defined for any tangent vector field $X$ on $S$ as
  $$ BX = D_XN~. $$
\item $\II$ is the second fundamental form of $S$, defined, for $X$ and $Y$ two vector fields on $S$, as
  $$ \II(X,Y)=I(BX,Y)=I(X,BY)~. $$
\item $\III$ is the third fundamental form of $S$, defined, for $X$ and $Y$ two vector fields on $S$, as
  $$ \III(X,Y)=I(BX,BY)~. $$
\end{itemize}

\begin{definition} \label{df:Isharp}
  Let $(u,\rho)$ be a locally convex equivariant immersion of a surface $S$ in $\HH^3$. We define the {\em metric at infinity} of $S$ as
  $$ I^\# = I((E+B)\cdot,(E+B)\cdot)~. $$
\end{definition}

A key property is that $I^\#$ is in the conformal class of the pull-back of the conformal class at infinity of $\HH^3$ on $S$ by the hyperbolic Gauss map, see e.g. \cite{volume}.

Note that different multiplicative coefficients can be found in different references. With the definition given here, the metric at infinity of a totally geodesic plane is hyperbolic, of constant curvature $-1$. The surface $S$ can be reconstructed from the metric at infinity by a construction due to Epstein \cite{epstein:envelopes}. (A parallel point of view is developed in \cite{horo}.)

The following lemma is similar to the analysis made in \cite{L5}.

\begin{lemma} \label{lm:Isharp}
  We denote by $D^\#$ the Levi-Civita connection of $I^\#$, by $K^\#$ its curvature, and by $J^\#$ its complex structure. We then have 
  \begin{equation}\label{eq:Dsharp}
    D^\#_xy=(E+B)^{-1}D_x((E+B)y)~,
  \end{equation}
  \begin{equation}\label{eq:Jsharp}   
    J^\#=(E+B)^{-1}J(E+B)~.
  \end{equation}
  \begin{equation}\label{eq:Ksharp}   
    K=\det(E+B)K^\#~.
  \end{equation}
  If $f:S\to \R$ is a smooth function,
  \begin{equation}
    \label{eq:gradient}
    D^\#f = (E+B)^{-2}Df~.
  \end{equation}
\end{lemma}

\begin{proof}
  The right-hand side of \eqref{eq:Dsharp} defines a connection, which is
  compatible with $I^\#$ since
  \begin{eqnarray*}
    x\cdot I^\#(y,z) & = & x\cdot I((E+B)y,(E+B)z) \\
    & = & I\big(D_x((E+B)y),(E+B)z\big)+I\big((E+B)y,D_x((E+B)z)\big) \\
    & = & I^\#(D^\#_xy,z)+I^\#(y,D^\#_xz)~.
  \end{eqnarray*}
  It is also torsion-free since $D$ is torsion-free and $B$ is Codazzi, so that
  \begin{eqnarray*}
      D^\#_xy-D^\#_yx-[x,y] & = & (E+B)^{-1}\big(D_x((E+B)y)-D_y((E+B)x)-(E+B)[x,y]\big) \\
                            & = & (E+B)^{-1}\,d^D(E+B)(x,y)=0~.
  \end{eqnarray*}
  So $D^\#$ is the Levi-Civita connection of $I^\#$.

  For \eqref{eq:Jsharp}, note that
  $$ \big((E+B)^{-1}J(E+B)\big)^2=(E+B)^{-1}J^2(E+B)=-E $$
  and
  \begin{eqnarray*}
    I^\#\big((E+B)^{-1}J(E+B)x,(E+B)^{-1}J(E+B)y\big) & = & I\big(J(E+B)x,J(E+B)y\big) \\
    &=&I\big((E+B)x,(E+B)y\big)=I^\#(x,y)~, 
  \end{eqnarray*}
  so $(E+B)^{-1}J(E+B)$ is a complex structure compatible with $I^\#$; it
  induces the same orientation as $J$ because $\det(E+B)>0$, hence it is
  $J^\#$.

  For \eqref{eq:Ksharp}, let $(e_1,e_2=Je_1)$ be a direct orthonormal frame for $I$ on an open subset $U\subset S$. By \eqref{eq:Jsharp}, $e^\#_1:=(E+B)^{-1}e_1$ and $e^\#_2:=(E+B)^{-1}e_2$ form a direct orthonormal frame for $I^\#$ on $U$, and \eqref{eq:Dsharp} gives
  $$ D^\#_xe^\#_1=(E+B)^{-1}D_xe_1~, $$
  so that the connection $1$-forms of $D$ and $D^\#$ in those frames are equal:
  $$ \omega^\#(x)=I^\#(D^\#_xe^\#_1,e^\#_2)=I(D_xe_1,e_2)=\omega(x)~. $$
  The structure equations $d\omega=-K\,da$ and $d\omega^\#=-K^\#\,da^\#$ therefore yield
  $K da=K^\# da^\#$, while
  $$ da^\#(e_1,e_2)=\det(E+B) da^\#(e^\#_1,e^\#_2)=\det(E+B)=\det(E+B) da(e_1,e_2)~, $$
  that is, $da^\#=\det(E+B) da$. Hence $K=\det(E+B)K^\#$ on $U$, and therefore on $S$.

  Finally, \eqref{eq:gradient} follows from the definition of the gradient and of $I^\#$. 
\end{proof}

\section{Infinitesimal rigidity}

We provide in this section the proofs of Lemmas \ref{lm:rigid} and \ref{lm:rigid*}, starting with the infinitesimal rigidity for the induced metric.

\subsection{Variations with prescribed $I$ and $c$} \label{ssc:variations}

We assume that $\Sigma\subset \HH^3$ is a locally strongly convex surface and use the usual notations $B, K, I, I^\#, K^\#$, etc. We denote by $D^\#$ the Levi-Civita connection of $I^\#$, and by $J^\#$ its complex structure.

We consider a first-order deformation of $\Sigma$, given by a first-order variation $\dot B$, such that $\dot I^\#$ is the sum of a conformal variation and a ``trivial'' variation associated to a vector field $v$. The variation of $I^\#$ is equal to
$$ \dot I^\#=I^\#(b\cdot, \cdot)+I^\#(\cdot, b\cdot)~, $$
where
$$ b=(E+B)^{-1}\dot B~, $$
and we assume that it can be written as the sum of a trivial deformation and of a conformal deformation, so as:
$$ \dot I^\#= \cL_vI^\# + 2\lambda I^\#~. $$
Since $\cL_vI^\#(x,y)=I^\#(D^\#_xv,y)+I^\#(x,D^\#_yv)$, the previous equation states that the self-adjoint parts of $b$ and of $D^\#v+\lambda E$ are equal (where $E$ is the identity). This can in turn be written as the existence of a function $\mu$ such that
\begin{equation}
  \label{eq:b}
  b = D^\#v + \lambda E + \mu J^\#~.
\end{equation}

\begin{lemma}
  $b$ satisfies:
  \begin{equation}
    \label{eq:codazzi}
    d^{D^\#}b=0~,
  \end{equation}
  \begin{equation}
    \label{eq:tr}
    \tr((E+B^{-1})b)=0~,
  \end{equation}
  \begin{equation}
    \label{eq:tr2}
    \tr((E+B)J^\#b)=0~.
  \end{equation}
\end{lemma}

\begin{proof}
  Since the metric $I$ is fixed along the deformation, so are its
  Levi-Civita connection $D$ and its curvature $K$. Differentiating at $t=0$
  the Codazzi equation $d^DB_t=0$, with $D$ constant, gives
  $$ d^D\dot B=0~. $$
  Differentiating the Gauss equation $\det B_t=1+K$ and using that $\frac{d}{dt}\det B_t=\det B_t\,\tr(B_t^{-1}\dot B_t)$ (with $\det B=1+K\neq0$) gives
  $$ \tr(B^{-1}\dot B)=0~. $$
  Finally each $B_t$ is $I$-self-adjoint and $I$ is fixed, so $\dot B$ is
  $I$-self-adjoint, and therefore
  $$ \tr(J\dot B)=0 $$
  by Lemma~\ref{lm:2d}(3).

  We now prove the three equations. For \eqref{eq:codazzi}, using \eqref{eq:Dsharp} and $(E+B)b=\dot B$,
  \begin{eqnarray*}
    d^{D^\#}b(x,y) & = & D^\#_x(by)-D^\#_y(bx)-b[x,y] \\
                   & =& (E+B)^{-1}\Big(D_x\big((E+B)by\big)-D_y\big((E+B)bx\big)-(E+B)b[x,y]\Big) \\
                   & = & (E+B)^{-1}\,d^D\dot B(x,y)~,
  \end{eqnarray*}
  which vanishes since $d^D\dot B=0$.

  For \eqref{eq:tr}, since $E+B^{-1}=B^{-1}(B+E)$ and $b=(E+B)^{-1}\dot B$,
  $$ \tr\big((E+B^{-1})b\big)=\tr\big(B^{-1}(E+B)(E+B)^{-1}\dot B\big)
  =\tr(B^{-1}\dot B)=0~. $$

  For \eqref{eq:tr2}, Equation \eqref{eq:Jsharp} gives
  $(E+B)J^\#=J(E+B)$, so that
  $$ (E+B)J^\#b=J(E+B)b=J\dot B~, $$
  and $\tr(J\dot B)=0$ because $\dot B$ is self-adjoint for $I$.
\end{proof}

\subsection{Expressing $v$ in terms of $\lambda$ and $\mu$}

The following lemma expresses $v$ in terms of $\lambda$ and $\mu$.

\begin{lemma}
  We have:
  \begin{equation}
     \label{eq:v2}
      K^\# v = D^\#\lambda +J^\#D^\#\mu~.
    \end{equation}

\end{lemma}

\begin{proof}
  We compute $d^{D^\#}$ of each term in \eqref{eq:b}.
  \begin{eqnarray*}
    (d^{D^\#}D^\#v)(x,y) & = & D^\#_xD^\#_yv-D^\#_yD^\#_xv-D^\#_{[x,y]}v \\
                       & = & R^\#(x,y)v \\
                       & = & (-K^\#J^\# v) da^\#(x,y)~, \\
    (d^{D^\#}\lambda E)(x,y) & = & D^\#_x(\lambda y) - D^\#_y(\lambda x) - \lambda [x,y] \\
                       & = & d\lambda(x)y-d\lambda(y)x \\
                       & = & J^\#(D^\#\lambda) da^\#(x,y)~, \\
    (d^{D^\#}\mu J^\#)(x,y) & = & D^\#_x(\mu J^\#y) - D^\#_y(\mu J^\# x) - \mu J^\#[x,y] \\
                       & = & d\mu(x)J^\#y - d\mu(y)J^\#x \\
    & = & -(D^\#\mu) da^\#(x,y)
  \end{eqnarray*}
  and the result follows from putting the three terms together with \eqref{eq:codazzi}.
\end{proof}

We can now express the first-order variation $\dot B$ in terms of $v, \lambda$ and $\mu$.

\begin{lemma}\label{lm:transl}
  Let $\dot B$, $v$, $\lambda$, $\mu$ be as above, and set $w=(E+B)v$. Then
  \begin{equation}
    \label{eq:dotBw}
    \dot B = Dw+\lambda(E+B)+\mu J(E+B)~,
  \end{equation}
  \begin{equation}
    \label{eq:Kw}
    Kw = \adj(E+B)\,D\lambda+(E+B)J\,D\mu~.
  \end{equation}
\end{lemma}


\begin{proof}
  For \eqref{eq:dotBw}, multiply \eqref{eq:b} on the left by $E+B$: since
  $\dot B=(E+B)b$, and using
  $$ D^\#_xv=(E+B)^{-1}D_x((E+B)v) $$
  and
  $$ J^\#=(E+B)^{-1}J(E+B)~, $$
  $$ \dot B \;=\; (E+B)\,D^\#v+\lambda (E+B)+\mu (E+B)J^\#= Dw+\lambda (E+B)+\mu J(E+B)~. $$

  For \eqref{eq:Kw}, note that $D^\#f=(E+B)^{-2}Df$, so \eqref{eq:v2} translates as
  \begin{eqnarray*}
    K^\#v &=& (E+B)^{-2}D\lambda+(E+B)^{-1}J(E+B)\cdot (E+B)^{-2}D\mu \\
          &=& (E+B)^{-2}D\lambda+(E+B)^{-1}J(E+B)^{-1}D\mu~.
  \end{eqnarray*}
  Multiplying by $E+B$ gives
  $$ K^\#w=(E+B)^{-1}D\lambda+J(E+B)^{-1}D\mu~. $$
  Multiplying by $\det(E+B)$, using $\det(E+B)\,(E+B)^{-1}=\adj(E+B)$ and \eqref{eq:Ksharp}, we obtain that
  $$ Kw=\adj(E+B)\,D\lambda+J\adj(E+B)\,D\mu~, $$
  while $J\adj(E+B)=(E+B)J$ by Lemma~\ref{lm:2d}(1).
\end{proof}

\subsection{An integral formula}

The following formula is the key to the proof of the main infinitesimal rigidity lemma.

\begin{proposition}\label{pr:identity}
  For every infinitesimal deformation of $(I,B)$ fixing the induced metric and
  the conformal class at infinity, with the notation of
  Lemma~\ref{lm:transl},
  \begin{equation}
    \label{eq:identity}
    \int_S K\, I(w,w)\, da \;=\; -2\int_S \det\dot B\, da
    \;+\; 2\int_S \det(E+B)\,(\lambda^2+\mu^2)\, da~.
  \end{equation}
\end{proposition}

No assumption on the sign of $K$ is made in this proposition.

\begin{proof}
  Recall that $E+B$ is self-adjoint, positive and Codazzi, and that $\dot B$
  is self-adjoint and Codazzi.

  Take the $I$-scalar product of \eqref{eq:Kw} with $w$ and integrate over $S$:
  $$ \int_S K\,I(w,w)\,da =\int_S I(\adj(E+B)D\lambda,w)\,da+\int_S I((E+B)JD\mu,w)\,da~. $$
  Since $\adj(E+B)$ is self-adjoint,
  $$ I(\adj(E+B)D\lambda,w)=d\lambda(\adj(E+B)w)~, $$
  so by Lemma~\ref{lm:div}(1) and (2),
  $$ \int_S I(\adj(E+B)D\lambda,w)\,da
  =-\int_S\lambda\,\dv(\adj(E+B)w)\,da
  =-\int_S\lambda\,\tr(\adj(E+B)\,Dw)\,da~. $$

  Similarly,
  $$ I((E+B)JD\mu,w)=I(JD\mu,(E+B)w)=-d\mu(J(E+B)w)~, $$
  so by Lemma~\ref{lm:div}(3),
  $$ \int_S I((E+B)JD\mu,w)\,da
  =\int_S\mu\,\dv(J(E+B)w)\,da
  =\int_S\mu\,\tr(J(E+B)\,Dw)\,da~. $$
  
  Now substitute
  $$ Dw=\dot B-\lambda(E+B)-\mu J(E+B) $$
  from \eqref{eq:dotBw}.
  Using Lemma~\ref{lm:2d},
  $$ \tr(\adj(E+B)(E+B))=2\det(E+B) $$ and
  $$ \tr(\adj(E+B)J(E+B))=0~, $$
  while by linearity of the adjugate and the linearized Gauss equation,
  \begin{equation}
    \label{eq:E+B}
    \tr(\adj(E+B)\dot B)=\tr((E+\adj(B))\dot B) =\tr\dot B+\det(B)\tr(B^{-1}\dot B)=\tr\dot B~,
  \end{equation}
  so
  $$ \tr(\adj(E+B)Dw)=\tr\dot B-2\lambda\det(E+B)~. $$

  Similarly,
  $$ \tr(J(E+B)\dot B)=\tr(J\dot B)+\tr(JB\dot B)=\tr(JB\dot B)~, $$
  while $\tr(J(E+B)^2)=0$ by Lemma~\ref{lm:2d}(3), and
  $\tr(J(E+B)J(E+B))=-2\det(E+B)$ by Lemma~\ref{lm:2d}(5), so
  $$ \tr(J(E+B)\,Dw)=\tr(JB\dot B)+2\mu\det(E+B)~. $$

  Putting the pieces together,
  \begin{equation}
    \label{eq:E1}
    \int_S K\,I(w,w)\,da
    =\int_S\big(-\lambda\,\tr\dot B+\mu\,\tr(JB\dot B)\big)\,da
    +2\int_S\det(E+B)\,(\lambda^2+\mu^2)\,da~.
  \end{equation}

  Since $\dot B$ is self-adjoint and Codazzi, Lemma~\ref{lm:div}(2) applied to $A=\dot B$ gives
  $$ \int_S\tr(\adj(\dot B) Dw)da=0~. $$
  Substituting \eqref{eq:dotBw} again:
  \begin{equation}
    \label{eq:interm}
    \int_S\tr(\adj(\dot B)(\dot B -\lambda(E+B)-\mu J(E+B)))da=0~. 
  \end{equation}
  Since $\adj(A)A=(\det A)E$ for any $A$ (Section \ref{ssc:conventions}),
  $$ \tr(\adj(\dot B)\dot B)=2\det\dot B~. $$
  By Lemma~\ref{lm:2d}(2) and \eqref{eq:E+B},
  $$ \tr(\adj(\dot B)(E+B))=\tr(\adj(E+B)\dot B)=\tr\dot B~. $$
  Finally, by Lemma~\ref{lm:2d}(1),(4),
  $$ \adj(\dot B)J(E+B)=(-J\dot BJ)J(E+B)=J\dot B(E+B)~, $$
  so
  $$\tr(\adj(\dot B)J(E+B)) = \tr(J\dot B(E+B))=\tr(J\dot B)+\tr(J\dot B\,B)=-\tr(JB\dot B)~. $$
  Replacing each term in \eqref{eq:interm},
  \begin{equation}
    \label{eq:E2}
    0=\int_S\big(2\det\dot B-\lambda\,\tr\dot B+\mu\,\tr(JB\dot B)\big)\,da~.
  \end{equation}
  Substituting \eqref{eq:E2} into \eqref{eq:E1} yields \eqref{eq:identity}.
\end{proof}

\subsection{Proof of Lemma \ref{lm:rigid}}

Consider a first-order deformation of $(E,u)$ preserving the induced metric on the image, given by $\dot B$, and let $w, \lambda$ and $\mu$ be as above. The left-hand side of \eqref{eq:identity} is non-positive since $K\leq 0$, while the right-hand side is non-negative since $\det(\dot B)\leq 0$ (because $\tr(B^{-1}\dot B)=0$). Therefore, both sides must be zero. This means that $\dot B=0$, and the first-order deformation is therefore trivial.

\subsection{Infinitesimal deformation preserving $\III$}

We now consider a first-order deformation preserving the third fundamental form $\III=I(B\cdot,B\cdot)$ rather than the induced metric. We denote by $D^*$, $K^*$, $J^*$, $da^*$ the Levi-Civita connection, curvature, complex structure and area form of $\III$, and we set
$$ B^*=B^{-1}~. $$

There is a close similarity between the proof of Lemma \ref{lm:rigid} and that of Lemma \ref{lm:rigid*}. The simplest way to understand this proximity is by formulating Lemma \ref{lm:rigid*} in terms of the surface $S^*$ dual to $S$ in the de Sitter space, see e.g. Sections \ref{ssc:dual} and \ref{sssc:dual}. This duality is however not necessary to follow or check the proof, and it is implicit below --- although the data marked with a $*$ can all be interpreted in terms of the dual surface.

\begin{lemma}\label{lm:dual}
  We have
  \begin{equation}\label{eq:dualdict}
    D^*_xy=B^{-1}D_x(By)~,\qquad J^*=B^{-1}JB~,\qquad
    da^*=\det(B)\,da~,\qquad K^*\det B=K~.
  \end{equation}
  Moreover $B^*$ is $\III$-self-adjoint, positive definite and Codazzi for
  $D^*$, the Gauss equation takes the form
  \begin{equation}\label{eq:dualgauss}
    \det B^*=1-K^*~,
  \end{equation}
  and
  \begin{equation}\label{eq:dualIsharp}
    I^\#=\III((E+B^*)\cdot,(E+B^*)\cdot)~,\qquad
    D^\#_xy=(E+B^*)^{-1}D^*_x\big((E+B^*)y\big)~,
  \end{equation}
  \begin{equation}\label{eq:dualJsharp}
    J^\#=(E+B^*)^{-1}J^*(E+B^*)~,\qquad
    K^\#\det(E+B^*)=K^*~.
  \end{equation}
\end{lemma}

\begin{proof}
  The identities \eqref{eq:dualdict} are proved exactly as
  \eqref{eq:Dsharp}, \eqref{eq:Jsharp} and \eqref{eq:Ksharp}, with $E+B$
  replaced by $B$, which is also $I$-self-adjoint, positive and Codazzi.

  $B^*$ is $\III$-self-adjoint because
  $\III(B^*x,y)=I(BB^*x,By)=I(x,By)$ and $\III(x,B^*y)=I(Bx,y)$ which are equal since $B$ is $I$-self-adjoint. It is positive definite because $\III(B^*x,x)=I(x,Bx)>0$ for $x\neq0$. It is Codazzi for $D^*$ since, by \eqref{eq:dualdict},
  $$ D^*_x(B^*y)=B^{-1}D_x(BB^*y)=B^{-1}D_xy~, $$
  so that
  $$ d^{D^*}B^*(x,y)=B^{-1}\big(D_xy-D_yx-[x,y]\big)=0~. $$
  For \eqref{eq:dualgauss}, the Gauss equation $\det B=1+K$ and the last
  identity of \eqref{eq:dualdict} give
  $$ \det B^*=\frac1{\det B}=\frac{\det B-K}{\det B}=1-K^*~. $$

  Finally, since $B(E+B^*)=E+B$, for any vector fields $x,y$ on $S$,
  $$ \III\big((E+B^*)x,(E+B^*)y\big)=I\big(B(E+B^*)x,B(E+B^*)y\big)
  =I\big((E+B)x,(E+B)y\big)=I^\#(x,y)~ $$
  As $E+B^*$ is $\III$-self-adjoint, positive and Codazzi for $D^*$, the remaining identities follow from Lemma~\ref{lm:Isharp} applied to the metric $\III$ and the operator $E+B^*$ in place of $I$ and $E+B$.
\end{proof}

We now consider a first-order deformation of $\Sigma$ for which $\III$ is fixed, given by the first-order variation $\dot B^*$ of $B^*$, and we assume as before that the variation of $I^\#$ is the sum of a trivial and of a conformal variation,
$$ \dot I^\#=\cL_vI^\#+2\lambda I^\#~. $$
By \eqref{eq:dualIsharp}, $\dot I^\#=I^\#(b^*\cdot,\cdot)+I^\#(\cdot,b^*\cdot)$
with
$$ b^*=(E+B^*)^{-1}\dot B^*~, $$
so, exactly as above, there is a function $\mu$ such that
\begin{equation}
  \label{eq:bstar}
  b^*=D^\#v+\lambda E+\mu J^\#~.
\end{equation}

\begin{lemma}\label{lm:bstareq}
  $b^*$ satisfies
  \begin{equation}\label{eq:dualcodazzi}
    d^{D^\#}b^*=0~,
  \end{equation}
  \begin{equation}\label{eq:dualtr}
    \tr\big((E+(B^*)^{-1})b^*\big)=0~,
  \end{equation}
  \begin{equation}\label{eq:dualtr2}
    \tr\big((E+B^*)J^\#b^*\big)=0~.
  \end{equation}
\end{lemma}

\begin{proof}
  Since $\III$ is fixed along the deformation, so are $D^*$ and $K^*$.
  Differentiating at $t=0$ the Codazzi equation $d^{D^*}B^*_t=0$, with $D^*$
  constant, gives
  $$ d^{D^*}\dot B^*=0~. $$
  Differentiating \eqref{eq:dualgauss}, that is $\det B^*_t=1-K^*$, and
  using the Jacobi formula together with $\det B^*=1-K^*\neq0$, gives
  $$ \tr\big((B^*)^{-1}\dot B^*\big)=0~. $$
  Finally each $B^*_t$ is $\III$-self-adjoint and $\III$ is fixed, so
  $\dot B^*$ is $\III$-self-adjoint and therefore
  $$ \tr(J^*\dot B^*)=0 $$
  by Lemma~\ref{lm:2d}(3).

  The three equations now follow as in the proof of the corresponding
  statement for $b$. For \eqref{eq:dualcodazzi}, using \eqref{eq:dualIsharp}
  and $(E+B^*)b^*=\dot B^*$,
  \begin{eqnarray*}
    d^{D^\#}b^*(x,y) & = & D^\#_x(b^*y)-D^\#_y(b^*x)-b^*[x,y] \\
    & = & (E+B^*)^{-1}\Big(D^*_x\big((E+B^*)b^*y\big)
      -D^*_y\big((E+B^*)b^*x\big)-(E+B^*)b^*[x,y]\Big) \\
    & = & (E+B^*)^{-1}\,d^{D^*}\dot B^*(x,y)=0~.
  \end{eqnarray*}
  For \eqref{eq:dualtr}, since $E+(B^*)^{-1}=(B^*)^{-1}(B^*+E)$,
  $$ \tr\big((E+(B^*)^{-1})b^*\big)
  =\tr\big((B^*)^{-1}(E+B^*)(E+B^*)^{-1}\dot B^*\big)
  =\tr\big((B^*)^{-1}\dot B^*\big)=0~. $$
  For \eqref{eq:dualtr2}, \eqref{eq:dualJsharp} gives
  $(E+B^*)J^\#=J^*(E+B^*)$, so that
  $(E+B^*)J^\#b^*=J^*(E+B^*)b^*=J^*\dot B^*$, whose trace vanishes.
\end{proof}

Equation \eqref{eq:v2} was obtained from \eqref{eq:b} and
\eqref{eq:codazzi} by a computation involving only the objects attached to
$I^\#$. Since \eqref{eq:bstar} and \eqref{eq:dualcodazzi} have exactly the
same form, the same computation gives
\begin{equation}\label{eq:dualv2}
  K^\#v=D^\#\lambda+J^\#D^\#\mu~.
\end{equation}

\begin{lemma}\label{lm:dualtransl}
  Let $\dot B^*$, $v$, $\lambda$, $\mu$ be as above, and set
  $w^*=(E+B^*)v$. Then
  \begin{equation}\label{eq:dualdotBw}
    \dot B^*=D^*w^*+\lambda(E+B^*)+\mu J^*(E+B^*)~,
  \end{equation}
  \begin{equation}\label{eq:dualKw}
    K^*w^*=\adj(E+B^*)\,D^*\lambda+(E+B^*)J^*\,D^*\mu~.
  \end{equation}
\end{lemma}

\begin{proof}
  For \eqref{eq:dualdotBw}, multiply \eqref{eq:bstar} on the left by
  $E+B^*$: since $\dot B^*=(E+B^*)b^*$, and using \eqref{eq:dualIsharp} and
  \eqref{eq:dualJsharp},
  $$ \dot B^*=(E+B^*)D^\#v+\lambda(E+B^*)+\mu(E+B^*)J^\#
  =D^*w^*+\lambda(E+B^*)+\mu J^*(E+B^*)~. $$

  For \eqref{eq:dualKw}, note that $D^\#f=(E+B^*)^{-2}D^*f$, so
  \eqref{eq:dualv2} translates as
  $$ K^\#v=(E+B^*)^{-2}D^*\lambda
  +(E+B^*)^{-1}J^*(E+B^*)^{-1}D^*\mu~. $$
  Multiplying by $E+B^*$ gives
  $K^\#w^*=(E+B^*)^{-1}D^*\lambda+J^*(E+B^*)^{-1}D^*\mu$, and multiplying by
  $\det(E+B^*)$, using $\det(E+B^*)(E+B^*)^{-1}=\adj(E+B^*)$ and
  \eqref{eq:dualJsharp},
  $$ K^*w^*=\adj(E+B^*)\,D^*\lambda+J^*\adj(E+B^*)\,D^*\mu~, $$
  while $J^*\adj(E+B^*)=(E+B^*)J^*$ by Lemma~\ref{lm:2d}(1).
\end{proof}

\subsection{A dual integral identity}

The following analog of Proposition \ref{pr:identity} is the key to the dual infinitesimal rigidity lemma. 

\begin{proposition}\label{pr:dualidentity}
  For every infinitesimal deformation of $(I,B)$ fixing the third
  fundamental form and the conformal class at infinity, with the notation of
  Lemma~\ref{lm:dualtransl},
  \begin{equation}\label{eq:dualidentity}
    \int_S K^*\,\III(w^*,w^*)\,da^*
    \;=\;-2\int_S\det\dot B^*\,da^*
    \;+\;2\int_S\det(E+B^*)\,(\lambda^2+\mu^2)\,da^*~.
  \end{equation}
\end{proposition}

No assumption on the sign of $K^*$ is made in this proposition.

\begin{proof}
  By Lemma~\ref{lm:dual}, $E+B^*$ is $\III$-self-adjoint, positive and
  Codazzi for $D^*$, and $\dot B^*$ is $\III$-self-adjoint and Codazzi for
  $D^*$. The proof of Proposition~\ref{pr:identity} uses only these
  properties of $E+B$ and $\dot B$, together with the two equations of
  Lemma~\ref{lm:transl}, the linearized Gauss equation, and
  Lemmas~\ref{lm:2d} and \ref{lm:div}; the latter two hold for any
  Riemannian metric on $S$. Since \eqref{eq:dualdotBw} and \eqref{eq:dualKw}
  have exactly the form of \eqref{eq:dotBw} and \eqref{eq:Kw}, with
  $(I,D,J,da,K,B,\dot B,w)$ replaced by
  $(\III,D^*,J^*,da^*,K^*,B^*,\dot B^*,w^*)$, the same computation applies
  verbatim. We indicate the two points where the linearized Gauss equation
  enters.

  In the first part of the proof, one needs
  $$ \tr\big(\adj(E+B^*)\dot B^*\big)=\tr\big((E+\adj(B^*))\dot B^*\big)=\tr\dot B^*+\det(B^*)\tr\big((B^*)^{-1}\dot B^*\big)=\tr\dot B^*~, $$
  by \eqref{eq:dualtr}; this replaces \eqref{eq:E+B}. In the second part, the same identity together with Lemma~\ref{lm:2d}(2) gives $\tr(\adj(\dot B^*)(E+B^*))=\tr\dot B^*$. With these substitutions, \eqref{eq:E1} is replaced by
  $$ \int_S K^* \III(w^*,w^*)da^*=\int_S\big(-\lambda\tr\dot B^*+\mu\tr(J^*B^*\dot B^*)\big)\,da^*+2\int_S\det(E+B^*)(\lambda^2+\mu^2) da^*~, $$
  and \eqref{eq:E2} by 
  $$ 0=\int_S\big(2\det\dot B^*-\lambda\,\tr\dot B^*
  +\mu\,\tr(J^*B^*\dot B^*)\big)\,da^*~, $$
  and substituting the second into the first yields \eqref{eq:dualidentity}.
\end{proof}

\subsection{Proof of Lemma \ref{lm:rigid*}}

Consider a first-order deformation of $(E,u)$ preserving the third fundamental form on the image and the conformal class at infinity, given by $\dot B^*$, and let $w^*,\lambda,\mu$ be as above. By \eqref{eq:dualdict}, $K^*\det B=K$ with $\det B>0$, so $K\leq0$ implies $K^*\leq0$ and the left-hand side of \eqref{eq:dualidentity} is non-positive. On the right-hand side, $\tr((B^*)^{-1}\dot B^*)=0$ and $B^*$ is positive definite, so $\det\dot B^*\leq0$, and $\det(E+B^*)>0$; both terms are therefore non-negative. Hence all three integrals vanish, so $\lambda\equiv\mu\equiv0$ and $\det\dot B^*\equiv0$, and $\dot B^*\equiv0$. Consequently $\dot B=-B\dot B^*B=0$.

\section{Properness}

\subsection{Properness for the induced metric}

We now move to the proof of Lemma \ref{lm:proper}. The proof has three steps. First (Section \ref{sssc:conv-ends}) we show that the pleated boundaries of the $E_n$ stay in a compact subset of $\cT_S$, and deduce, using Thurston's parametrisation of hyperbolic ends and a theorem of Dumas--Wolf, that the ends $E_n$ themselves converge after extraction. Second (Section \ref{sssc:labourie}) we recall Labourie's a priori estimates for convex surfaces in a hyperbolic end \cite{L4} and his compactness theorem for elliptic isometric immersions \cite{L1}. Finally (Section \ref{sssc:proof}) we combine them.

\subsubsection{Convergence of the ends}
\label{sssc:conv-ends}

We will use here the notations introduced at the beginning of Section \ref{ssc:ends}.

\begin{theorem}[Thurston; see \cite{kamishima-tan}]\label{tm:thurston}
The grafting map $\mathrm{Gr}\colon\cML_S\times\cT_S\to\cP_S$ is a homeomorphism.
A hyperbolic end $E$ is determined by $\sigma(E)$, and $\sigma(E)=\mathrm{Gr}_{l(E)}(m(E))$. Therefore $c(E)=\mathrm{gr}_{l(E)}(m(E))$, where $\mathrm{gr}=\pi_0\circ\mathrm{Gr}$ is conformal grafting, and the map sending $E\in\cE_S$ to $c(E)\in \cT_S$ is continuous.
\end{theorem}

\begin{definition}\label{df:convergence}
A sequence of marked hyperbolic ends $(E_n)_{n\in \N}$ \emph{converges} to $E$ if $(m(E_n),l(E_n))\to(m(E),l(E))$ in $\cT_S\times\cML_S$, or equivalently, by Theorem \ref{tm:thurston}, if $\sigma(E_n)\to\sigma(E)$ in $\cP_S$. 
\end{definition}

Note that here by ``hyperbolic end'' we always consider hyperbolic metrics on $S\times (0,\infty)$ which are isometric to a hyperbolic end.

Convergence in this sense implies geometric convergence: the pleated boundaries converge as equivariant maps $\tilde S\to\HH^3$ up to conjugation, so the metric completions converge uniformly on $R$-neighborhood of the $\partial_0E_n$, for all $R>0$.

\begin{theorem}[Dumas--Wolf \cite{dumas-wolf}]\label{tm:dw}
For every $Y\in\cT_S$, the map $\lambda\mapsto\mathrm{gr}_\lambda(Y)$ is a homeomorphism from $\cML_S$ to $\cT_S$.
\end{theorem}

\begin{corollary}\label{cr:joint}
The map $G\colon\cT_S\times\cML_S\to\cT_S\times\cT_S$, defined by $G(m,l)=(m,\mathrm{gr}_l(m))$, is a homeomorphism. If $(m_n,l_n)_{n\in \N}$ is a sequence with $m_n$ in a compact subset of $\cT_S$ and $\mathrm{gr}_{l_n}(m_n)$ in a compact subset of $\cT_S$, then $(m_n,l_n)_{n\in \N}$ has a convergent subsequence.
\end{corollary}

\begin{proof}
$G$ is continuous by Theorem~\ref{tm:thurston}, and bijective by Theorem \ref{tm:dw} applied fibrewise. Both $\cT_S\times\cML_S$ and $\cT_S\times\cT_S$ are homeomorphic to $\R^{12g-12}$, so $G$ is a homeomorphism by invariance of domain, hence proper.
\end{proof}

\begin{lemma}\label{lm:retraction}
Let $E$ be a hyperbolic end and $r_E\colon E\to\partial_0 E$ the nearest-point retraction onto the pleated boundary. Then $r_E$ is $1$-Lipschitz for the path metrics. Consequently, if $u\colon (S,h)\to E$ is a good isometric embedding, then for every closed curve $\gamma$ in $S$,
$$ \ell_{m(E)}(\gamma)\le\ell_h(\gamma)~. $$
\end{lemma}

\begin{proof}
  Suppose first that the measured lamination $l(E)$ is rational, that is, its support is a disjoint union of simple closed curves. In this case, $E$ is constructed by an explicit gluing of hyperbolic pieces associated to ``flat'' parts of $(S,m(E))$ which are the connected components of the complement of the support of $l(E)$, and to the closed geodesics which are the closed curves in the support of $l(E)$. The fact that $r_E$ is contracting follows from the explicit form of the hyperbolic metric on those ``pieces''. The general case follows from the density of rational measured laminations in $\cML_S$.
\end{proof}



\begin{proposition}\label{pr:endsconverge}
Let $(h_n)_{n\in \N}$ be metrics on $S$ converging smoothly to $h$, and let $(E_n,u_n)\in\cE$ with $u_n:(S,h_n)\to E_n$ good isometric embeddings. Then the $m(E_n)$ lie in a compact subset of $\cT_S$. If moreover $c(E_n)\to c$ in $\cT_S$, then, after extracting a subsequence, $(E_n)_{n\in \N}$ converges to a hyperbolic end $E$ with $c(E)=c$.
\end{proposition}

\begin{proof}
Fix a filling system $\alpha_1,\dots,\alpha_k$. Since $h_n\to h$, $L:=\sup_{n\in \N,1\leq i\leq k}\ell_{h_n}(\alpha_i)<\infty$, and Lemma \ref{lm:retraction} shows that $\ell_{m(E_n)}(\alpha_i)\le L$. This shows the first claim since the set of points of $\cT_S$ for which the $\alpha_i, 1\leq i\leq k$, have bounded hyperbolic length is compact.

If $c(E_n)\to c$, then $\mathrm{gr}_{l(E_n)}(m(E_n))=c(E_n)$ stays in a compact set, and Corollary \ref{cr:joint} shows that $(m(E_n),l(E_n))$ subconverges to some $(m,l)$. Let $E$ be the end with $\sigma(E)=\mathrm{Gr}_l(m)$. Then $E_n\to E$ and, by continuity of grafting, $c(E)=\lim c(E_n)=c$.
\end{proof}


\subsubsection{Labourie's a priori estimates and compactness theorem}
\label{sssc:labourie}

The following is \cite[Prop.~4.1]{L4}. There, a hyperbolic end is called a geometrically finite end, and $S_k\subset E$ denotes the unique incompressible surface of constant curvature $k\in(-1,0)$ in $E$ \cite{L6}, so that the $S_k$ foliate $E$. In \cite[Prop.~4.1]{L4} the surfaces need to be ``incompressible'', but the good surfaces considered here satisfy this hypothesis.

\begin{theorem}[Labourie]\label{tm:labourie41}
Let $E$ be a hyperbolic end and $\Sigma\subset E$ an embedded, locally convex, incompressible surface, with induced metric $g$ of curvature $K\ge k_0>-1$. Then:
\begin{enumerate}[(i)]
\item $\Sigma$ lies outside $S_{k_0}$,
\item if $\Sigma$ lies outside $S_k$ then $k\le 2\pi\chi(S)/\mathrm{Area}(g)$, and if $S_k$ lies outside $\Sigma$ then $k\ge 2\pi\chi(S)/\mathrm{Area}(g)$,
\item $C_1\le d(x,\partial_0 E)\le C_2$ for all $x\in\Sigma$;
\item $\int_\Sigma H\,da\le C_3$.
\end{enumerate}
Here $C_1=d(S_{k_0},\partial_0 E)$, $C_2=\mathrm{diam}(g)+\sup_{S_{k_1}}d(\cdot,\partial_0 E)$ for any $k_1\in(2\pi\chi(S)/\mathrm{Area}(g),0)$, and $C_3$ depends only on $C_1,C_2,k_0,\mathrm{Area}(g)$ and $E$.
\end{theorem}


  We first make precise the geometric convergence mentioned after Definition \ref{df:convergence}, in the form in which it will be used.

\begin{lemma}\label{lm:geomconv}
Let $(E_n)_{n\in \N}$ converge to $E$ and let $N\subset E$ be a compact domain with smooth boundary. For $n$ large there are smooth embeddings $j_n:N\to E_n$ such that $j_n^*g_{E_n}\to g_E$ in $C^\infty(N)$ and $d_{E_n}(j_n(x),\partial_0E_n)\to d_E(x,\partial_0E)$ uniformly on $N$.
\end{lemma}

\begin{proof}[Sketch of proof]
Since $(m(E_n),l(E_n))\to(m(E),l(E))$, the lifts to $\HH^3$ of the pleated boundaries converge as $\rho_n$-equivariant maps $\tilde S\to\HH^3$, uniformly on compact sets, with $\rho_n\to\rho$. The normal flow from $\partial_0 E_n$ therefore provides homeomorphisms $\iota_n$ from a neighborhood of $\bar N$ in $E$ onto open subsets of $E_n$ such that the developing maps satisfy $\mathrm{dev}_{E_n}\circ\tilde\iota_n\to\mathrm{dev}_E$ uniformly on the preimage of $N$ in $\tilde E$, and such that $d_{E_n}(\iota_n(x),\partial_0E_n)\to d_E(x,\partial_0E)$ uniformly on $N$. Since $\mathrm{dev}_{E_n}$ is injective on balls of radius $\frac12 d(\cdot,\partial_0\tilde E_n)$, the local inverses $(\mathrm{dev}_{E_n})^{-1}\circ\mathrm{dev}_E$ patch to a local isometry $\tilde j_n:\tilde N\to\tilde E_n$ which is $C^0$-close to $\tilde\iota_n$ and almost equivariant, in the sense that $\tilde j_n\circ\gamma$ and $\gamma\circ\tilde j_n$ differ by isometries of $\HH^3$ close to the identity, for $\gamma$ in a finite generating set of $\pi_1S$. Averaging with a partition of unity yields an equivariant smooth map $\tilde j_n$ which is $C^\infty$-close to a local isometry, hence descends to the required $j_n$.
\end{proof}

Below, $\mathrm{Id}$ denotes the identity endomorphism (denoted $E$ in Section \ref{ssc:conventions}), to avoid confusion with the end.

\begin{lemma}\label{lm:kuniform}
Let $k\in(-1,0)$.
\begin{enumerate}[(i)]
\item For every hyperbolic end $E$ and every good surface $\Sigma\subset E$ with curvature $K\leq k$,
$$ \sup_{x\in\Sigma}d(x,\partial_0E)\;\leq\;\operatorname{artanh}\sqrt{1+k}~. $$
In particular $\sup_{S_k(E)}d(\cdot,\partial_0E)\leq \operatorname{artanh}\sqrt{1+k}$ for every $E$.
\item If $E_n\to E$, then $S_k(E_n)\to S_k(E)$: with $N$ a compact neighborhood of $S_k(E)$ in $E$ and $j_n$ as in Lemma \ref{lm:geomconv}, for $n$ large $S_k(E_n)=j_n(\Sigma_n)$ for closed surfaces $\Sigma_n\subset N$ converging to $S_k(E)$ in $C^{2,\alpha}$. In particular $d(S_k(E_n),\partial_0E_n)\to d(S_k(E),\partial_0E)$.
\end{enumerate}
\end{lemma}

\begin{proof}
(1) Let $x_1\in\Sigma$ be a point where $d(\cdot,\partial_0E)$ is maximal on $\Sigma$, let $t_1$ be this maximum, and let $y_1\in\partial_0E$ be a nearest point to $x_1$. Work in $\tilde E$ and in $\HH^3$ through the developing map, near the segment $[y_1,x_1]$. Let $P$ be the support plane at $y_1$ of the (locally convex) developed pleated surface which is orthogonal to $[y_1,x_1]$, let $H_P$ be the closed half-space bounded by $P$ containing the pleated surface near $y_1$, and let $P_{t_1}$ be the surface at distance $t_1$ from $P$ on the side of $x_1$, so that $x_1\in P_{t_1}$. For $z\in\Sigma$ near $x_1$, a nearest point $y_z\in\partial_0E$ is near $y_1$, hence $y_z\in H_P$, and $d(z,P)\leq d(z,y_z)\leq t_1$. Thus $\Sigma$ lies, near $x_1$, in the closed $t_1$-neighborhood of $H_P$, and touches its boundary $P_{t_1}$ at $x_1$. Since $\Sigma$ is convex and $P_{t_1}$ is the totally umbilical surface with $B=\tanh(t_1)\mathrm{Id}$, the maximum principle gives $B_\Sigma(x_1)\geq\tanh(t_1)\mathrm{Id}$, and therefore
$$ 1+k\;\geq\;1+K(x_1)\;=\;\det B_\Sigma(x_1)\;\geq\;\tanh^2(t_1)~. $$

(2) Let $u_0:S\to E$ be the embedding with image $S_k(E)$, let $\nu$ be its unit normal pointing away from $\partial_0E$, and let $B_0=\nabla\nu$ be its shape operator, so that $B_0>0$ and $\det B_0=1+k$. For a hyperbolic metric $g$ on $N$ close to $g_E$ in $C^\infty$ and $\psi\in C^{2,\alpha}(S)$ small, let $u_{g,\psi}(x)=\exp^g_{u_0(x)}(\psi(x)\nu(x))$, and set
$$ F(g,\psi)\;=\;\det B_g(u_{g,\psi})-(1+k)\;\in\;C^{0,\alpha}(S)~, $$
where $B_g(u_{g,\psi})$ is the shape operator of $u_{g,\psi}$ for $g$. Then $F(g_E,0)=0$, $F$ is continuous in $g$ and $C^1$ in $\psi$, with $D_\psi F$ continuous in $(g,\psi)$. At $(g_E,0)$ the first variation of the shape operator under the normal variation $\psi\nu$ is $-\mathrm{Hess}\,\psi+\psi(\mathrm{Id}-B_0^2)$, so that, using $\tr(B_0^{-1})=\tr(B_0)/\det B_0$,
$$ L\psi\;:=\;D_\psi F(g_E,0)\psi\;=\;\det B_0\,\tr\big(B_0^{-1}\dot B\big)\;=\;-(1+k)\,\tr\big(B_0^{-1}\mathrm{Hess}\,\psi\big)\;-\;k\,H_0\,\psi~. $$
$L$ is elliptic since $B_0^{-1}>0$, and its zero order coefficient $-kH_0$ is positive. If $L\psi=0$ and $\psi$ has a positive maximum at $x$, then $\mathrm{Hess}\,\psi(x)\leq 0$, so $-(1+k)\tr(B_0^{-1}\mathrm{Hess}\,\psi)(x)\geq 0$ and $L\psi(x)>0$, a contradiction; similarly at a negative minimum. So $L$ is injective and, being a second order elliptic operator on a closed surface (of index zero), it is an isomorphism $C^{2,\alpha}(S)\to C^{0,\alpha}(S)$. By the implicit function theorem with parameters, there are a neighborhood $\mathcal G$ of $g_E$ in $C^\infty(N)$ and a continuous map $g\mapsto\psi_g\in C^{2,\alpha}(S)$ on $\mathcal G$, with $\psi_{g_E}=0$ and $F(g,\psi_g)=0$.

For $n$ large, $g_n:=j_n^*g_{E_n}\in\mathcal G$ by Lemma \ref{lm:geomconv}. The surface $\Sigma_n:=u_{g_n,\psi_{g_n}}(S)$ is embedded and isotopic to $u_0(S)$, its shape operator is positive definite (an open condition, satisfied at $(g_E,0)$), and it has curvature $-1+\det B=k$ for $g_n$. Hence $j_n(\Sigma_n)$ is a closed embedded locally convex incompressible surface of constant curvature $k$ in $E_n$, and $j_n(\Sigma_n)=S_k(E_n)$ by uniqueness \cite{L6}. Since $\psi_{g_n}\to 0$ in $C^{2,\alpha}$, $\Sigma_n\to S_k(E)$ in $C^{2,\alpha}$, and the last assertion follows from the uniform convergence of $d_{E_n}(j_n(\cdot),\partial_0E_n)$ to $d_E(\cdot,\partial_0E)$ on $N$.
\end{proof}

\begin{remark}\label{rk:uniform}
  For a sequence $E_n\to E$ (Definition \ref{df:convergence}) and metrics $h_n\to h$, one needs the constants $C_1,C_2,C_3$ of Theorem \ref{tm:labourie41} to be uniform in $n$. Fix $k_0\in(-1,0)$ with
  $$ k_0\leq\inf_n\inf_S K_{h_n} $$
  and $ k_1\in(-1,0)$ with
  $$ k_1>\sup_n 2\pi\chi(S)/\area(h_n)~. $$
  (This is possible since $h_n\to h$ and $\inf K_h\in(-1,0)$ by the Gauss--Bonnet Theorem.)

  Then $d(S_{k_0}(E_n),\partial_0E_n)$ converges to $d(S_{k_0}(E),\partial_0E)>0$ by Lemma \ref{lm:kuniform}(ii), hence is bounded below, while $\mathrm{diam}(h_n)+\sup_{S_{k_1}(E_n)}d(\cdot,\partial_0E_n)\leq \mathrm{diam}(h_n)+\operatorname{artanh}\sqrt{1+k_1}$ by Lemma \ref{lm:kuniform}(i), hence is bounded above, so that $C_1$ and $C_2$ can be chosen to be uniform constants. Given these, the injectivity radius of $E_n$ is bounded below on $\{C_1\leq d(\cdot,\partial_0E_n)\leq C_2\}$ uniformly in $n$, and $C_3$ is uniform by the argument of \cite[4.1(iv)]{L4} (cover $\Sigma$ by a bounded number of $\varepsilon$-balls and apply \cite[Lemme~5.4(iii)]{L1} in each).
\end{remark} 

We will use below the following result.

\begin{theorem}[Labourie {\cite[Th\'eor\`eme~D]{L1}}]\label{tm:labourieD}
  Let $(M,g)$ be a Riemannian $3$-manifold, $\Sigma$ a surface, and let $(f_n)_{n\in \N}$ be a sequence of isometric immersions $f_n:(\Sigma,g_n)\to(M,g)$ with $\det B_n\ge\varepsilon>0$, such that $g_n\to g_0$ in $C^\infty_{\rm loc}$ and $f_n\to f_0$ in $C^0_{\rm loc}$.
  If $\int_O H_n\,da_n$ is bounded for a neighbourhood $O$ of $x\in\Sigma$, a subsequence of $(f_n)$ converges in $C^\infty$ on a neighbourhood of $x$.
\end{theorem}

\begin{proposition}\label{pr:endsimplyconv}
Assume that $h_n\to h$ smoothly, with $K_h>-1$, and let $(E_n,u_n)_{n\in \N}$ be a sequence of elements of $\cE$ with $u_n:(S,h_n)\to E_n$ good isometric embeddings, and $(E_n)_{n\in \N}$ converging to $E$. Then, after extraction of a subsequence,  $(E_n,u_n)_{n\in \N}$ converges smoothly to $(E,u)\in\cE$, with $u:(S,h)\to E$ a good isometric embedding.
\end{proposition}

\begin{proof}
    By Lemma~\ref{lm:equiv}, lift $(E_n,u_n)$ to equivariant immersions $(\phi_n,\rho_n)$ of $\tilde S$ into $\HH^3$. We can normalise the $\phi_n$ by composing on the left by a global isometry so that $j^1\phi_n(x_0)$ converges for a fixed $x_0\in\tilde S$. For each $\gamma$ in a finite generating set of $\pi_1S$, $\rho_n(\gamma)$ sends the $1$-jet of $\phi_n$ at $x_0$ to its $1$-jet at $\gamma x_0$, and $d(\phi_n(x_0),\phi_n(\gamma x_0))\leq d_{\tilde h_n}(x_0,\gamma x_0)$ is bounded, so after extracting a subsequence $\rho_n(\gamma)$ converges to an isometry $\rho'(\gamma)$; this defines a representation $\rho':\pi_1S\to\PSL(2,\C)$.

  The $\phi_n$ are $1$-Lipschitz for $\tilde h_n\to\tilde h$, so the Ascoli Theorem shows that, after extracting a sub-sequence, $\phi_n\to \phi$ in $C^0_{\rm loc}$. The immersions are uniformly elliptic, $\det B_n=1+K_{h_n}\ge\varepsilon>0$, and Theorem \ref{tm:labourie41}(iv) bounds $\int_S H_n\,da_n$, hence Theorem \ref{tm:labourieD} applies at every point and a diagonal extraction gives that $\phi_n\to\phi$ in $C^\infty_{\rm loc}(\tilde S)$.

    The limit $\phi$ is an isometric immersion of $(\tilde S,\tilde h)$ with $\det B=1+K_h>0$, hence locally strongly convex, and passing to the limit in $\phi_n\circ\gamma=\rho_n(\gamma)\circ\phi_n$ shows that it is $\rho'$-equivariant. By Lemma~\ref{lm:equiv} it defines a good embedding $u\colon S\to E'$ in a hyperbolic end $E'$. The $\bCP^1$-structure of $E'$ is that of $(\phi,\rho')$, i.e. the limit of the $\bCP^1$-structures of $(\phi_n,\rho_n)$, which are $\sigma(E_n)\to\sigma(E)$, and by Theorem~\ref{tm:thurston}, $E'=E$.
\end{proof}

\subsubsection{Proof of Lemma \ref{lm:proper}}
\label{sssc:proof}

Let $h_n\to h$, $(E_n,u_n)$ and $c_n\to c$ be as in Lemma~\ref{lm:proper}.
By Proposition~\ref{pr:endsconverge}, after extraction $E_n\to E$ with
$c(E)=c$. By Proposition~\ref{pr:endsimplyconv}, after a further extraction
$(E_n,u_n)\to(E,u)$ smoothly, with $u\colon(S,h)\to E$ a good isometric
embedding. Since $c(E)=c$, this is the conclusion of Lemma~\ref{lm:proper}.
\qed

Note that the only places where $K>-1$ is used are Theorem~\ref{tm:labourie41} (through $k_0>-1$) and Theorem~\ref{tm:labourieD} (uniform ellipticity), but the hypothesis $K\le0$ of Theorem~\ref{tm:main} plays no role in properness.

\subsection{Properness for $\III$}

The proof of Lemma~\ref{lm:proper*} follows the same three steps as the proof of Lemma~\ref{lm:proper} in Section \ref{sssc:conv-ends}--\ref{sssc:proof}, with the hyperbolic end replaced, in the first step, by the de Sitter spacetime dual to it. After recalling in Section \ref{sssc:dual} the duality between locally convex surfaces in a hyperbolic end and spacelike Cauchy surfaces in the dual de Sitter spacetime, we show in Section \ref{sssc:ends*} that the bending laminations of the $E_n$ stay in a compact subset of $\cML_S$, and deduce, using Scannell--Wolf \cite{scannell-wolf} in place of Dumas--Wolf, that $(E_n)_{n\in \N}$ converges after extracting a subsequence. This corresponds directly to Section \ref{sssc:conv-ends} under duality: the pleated boundary is replaced by the initial singularity (the $\R$-tree dual to the bending lamination), and the $1$-Lipschitz retraction bounds the intersection numbers $i(l(E_n),\cdot)$ instead of the lengths $\ell_{m(E_n)}(\cdot)$. The second step, in Section \ref{sssc:dualest}, applies Labourie's a priori estimates (Theorem~\ref{tm:labourie41}) to the hyperbolic surfaces $u_n(S)$ themselves. The only point requiring care is the upper bound on the distance to the pleated boundary, which in Theorem~\ref{tm:labourie41}(iii) is obtained from the diameter of the induced metric, a quantity not controlled by the third fundamental form; it is here that the hypothesis $K^*<0$ enters, through the maximum principle of Lemma~\ref{lm:kuniform}(i). The third step, in Section \ref{sssc:compact*}, uses the Lorentzian analogue of Labourie's compactness theorem, proved in \cite{these}, on the dual surfaces, since it is the third fundamental forms, and not the induced metrics, that are known to converge.

\subsubsection{Dual surfaces and the de Sitter spacetime of an end}
\label{sssc:dual}

We use the following facts; see \cite{BBZ2,mess,scannell,benedetti-bonsante} and \cite{fillastre-seppi:spherical}, and Lemma \ref{lm:dual} for $B^*=B^{-1}$, $\det B^*=1-K^*$ and $I^\#=\III((E+B^*)\cdot,(E+B^*)\cdot)$.

\begin{enumerate}
\item[(D1)] To a hyperbolic end $E$ is associated a maximal globally hyperbolic, future complete, spatially compact de Sitter spacetime $E^*$, homeomorphic to $S\times(0,\infty)$, whose future conformal boundary carries the $\bCP^1$-structure $\sigma(E)$. Every such spacetime corresponds in this manner to a unique hyperbolic end.
\item[(D2)] The cosmological time $T:E^*\to(0,\infty)$ is $C^{1,1}$, its level sets $B^*_t$ are the duals of the equidistant surfaces $B_t$ from $\partial_0 E$ in $E$, and the past-directed gradient lines of $T$ retract $E^*$ onto the initial singularity, which (in the universal cover) is the $\R$-tree $T_{l(E)}$ dual to the bending lamination, with the metric given by the transverse measure. When $l(E)$ is rational, the pieces of $B^*_t$ dual to totally geodesic pieces of $\partial_0 E$ have induced metric $\sinh^2(t) m(E)$, while on the strip dual to a bending line of length $\ell$ and angle $\theta$ it is $\sinh^2(t) ds^2+\cosh^2(t)d\alpha^2$, $(s,\alpha)\in[0,\ell]\times[0,\theta]$.
\item[(D3)] If $(E,u)\in\cE$, the unit normal $u^*=N\colon S\to E^*$ is a spacelike embedding whose image is a Cauchy surface, with induced metric $I^*=\III$ and shape operator $B^*=B^{-1}$. Conversely the unit timelike normal of an equivariant uniformly locally strongly convex spacelike immersion of $\tilde S$ in $dS^3$ is an equivariant locally strongly convex immersion in $\HH^3$ whose third fundamental form is the induced metric of the former (see e.g. \cite{BBZ2}).
\item[(D4)] The translation length of $g\in\pi_1S$ acting on $T_l$ is $i(l,g)$ \cite{morgan-shalen}.
\item[(D5)] Pointwise,
  $$ H^* da^*=\tr(B^{-1})\det B da=\tr(\adj B) da=H da~. $$
\end{enumerate}

\subsubsection{Convergence of ends}
\label{sssc:ends*}

We prove in this section that the sequence of hyperbolic ends associated to a sequence of third fundamental forms and conformal structures at infinity converges (after taking a subsequence), see Proposition \ref{pr:endsconverge*} below.

The following lemma is closely related to similar estimates, see e.g. \cite{belraouti2}, \cite[\S 6]{cyclic}. We include a proof for completeness.

\begin{lemma}\label{lm:retraction*}
Let $E$ be a hyperbolic end and let $r^*:\widetilde E^*\to T_{l(E)}$ be the retraction of (D2). For every spacelike curve $c$ in $\widetilde E^*$,
$$ \ell_{T_{l(E)}}(r^*\circ c)\leq \ell(c)~. $$
Consequently, if $(E,u)\in\cE$ has third fundamental form $h^*$, then for every closed curve $g$ on $S$,
$$ i(l(E),g)\leq\ell_{h^*}(g)~. $$
\end{lemma}

\begin{proof}
Write the metric of $\widetilde E^*$ as $-dT^2+g_T$ in the $C^{1,1}$ coordinates given by $T$ and its gradient lines. Two properties are needed below: that $g_T$ is non-decreasing in $T$, and that the collapsing map $B^*_{T_0}\to T_{l(E)}$ is $1$-Lipschitz for $g_{T_0}$, for every $T_0>0$. When $l(E)$ is rational both follow from the explicit description in (D2): $\sinh^2(t)$ and $\cosh^2(t)$ are increasing, the map to the tree collapses the pieces with metric $\sinh^2(T_0)m(E)$ to vertices, and sends each strip to its edge by $(s,\alpha)\mapsto \alpha$, which has norm $1/\cosh(T_0)\leq 1$. The same two properties still hold when $l(E)$ is not rational by the density of rational measured laminations.

Let $c:s\to (x(s),T(s))$, $s\in \R/L\Z$, be a spacelike parameterized closed curve, and let $T_0=\min_cT$. Its length is equal to 
$$ \ell(c)=\int_0^L\sqrt{g_T(x'(s),x'(s))-T'(s)^2}ds \leq \int_0^L|x'(s)|_{g_T}ds~, $$
and projecting $c$ along the gradient lines onto $B^*_{T_0}$ gives a curve of length
$$ \int_0^L|x'(s)|_{g_{T_0}}ds\leq \int_0^L|x'(s)|_{g_T}ds~. $$
Composing with the collapsing map shows that $r^*$ does not increase the length of spacelike curves.

For the second claim, $r^*\circ u^*$ is $\rho(E)$-equivariant and $1$-Lipschitz from $(\tilde S,\tilde h^*)$ to $T_{l(E)}$. Applying it to a point $x$ on the axis of the isometry corresponding to $g$ and using (D4) gives, for any closed geodesic $g$ on $S$ corresponding to an element $\gamma$ of $\pi_1S$,
$$ i(l(E),g)\leq d_{T_{l(E)}}(r^*u^*(x),\gamma\cdot r^*u^*(x))\leq d_{\tilde h^*}(x,\gamma x)=\ell_{h^*}(g)~. $$
\end{proof}

Lemma \ref{lm:retraction*} is the exact analogue of Lemma \ref{lm:retraction}, with the lengths for $m(E)$ replaced by the intersection numbers with $l(E)$. The compactness statement replacing the compactness of the set of hyperbolic metrics of bounded length on a filling system is the following.

\begin{lemma}\label{lm:compactML}
Let $\alpha_1,\dots,\alpha_k$ be simple closed curves filling $S$. For every $L>0$, the set $\{l\in\cML_S\tq i(l,\alpha_i)\leq L,\ i=1,\dots,k\}$ is compact.
\end{lemma}

\begin{proof}
The function $l\mapsto\max_ii(l,\alpha_i)$ is continuous, homogeneous of degree one, and positive on $\cML_S\setminus\{0\}$ because the $\alpha_i$ fill $S$. A sublevel set of such a function is compact.
\end{proof}

\begin{theorem}[Scannell--Wolf \cite{scannell-wolf}]\label{tm:sw}
For every $l\in\cML_S$, the map $m\mapsto \mathrm{gr}_l(m)$ is a homeomorphism from $\cT_S$ to $\cT_S$.
\end{theorem}

\begin{corollary}\label{cr:joint*}
The map $G^*:\cML_S\times\cT_S\to\cML_S\times\cT_S$, defined by $G^*(l,m)=(l,\mathrm{gr}_l(m))$, is a homeomorphism. 
\end{corollary}

\begin{proof}
$G^*$ is continuous by Theorem \ref{tm:thurston}, and bijective by Theorem \ref{tm:sw} applied fibrewise. Both $\cML_S\times\cT_S$ and its image are homeomorphic to $\R^{12g-12}$, so $G^*$ is a homeomorphism by invariance of domain, hence proper.
\end{proof}

Therefore, if $(l_n,m_n)_{n\in \N}$ is a sequence with $l_n$ in a compact subset of $\cML_S$ and $\mathrm{gr}_{l_n}(m_n)$ in a compact subset of $\cT_S$, then $(l_n,m_n)_{n\in \N}$ has a convergent subsequence.

\begin{proposition}\label{pr:endsconverge*}
Let $(h^*_n)_{n\in \N}$ be metrics on $S$ converging smoothly to $h^*$, and let $(E_n,u^*_n)_{n\in \N}$ be a sequence of elements of $\cE$, where $u^*_n:S\to E_n$ is a good embedding with third fundamental form $h^*_n$. Then the $l(E_n)$ lie in a compact subset of $\cML_S$. If moreover $c(E_n)\to c$ in $\cT_S$, then, after extracting a subsequence, $(E_n)_{n\in \N}$ converges to a hyperbolic end $E$ with $c(E)=c$.
\end{proposition}

\begin{proof}
Fix a filling system $\alpha_1,\dots,\alpha_k$. Since $h^*_n\to h^*$, $L:=\sup_{n,i}\ell_{h^*_n}(\alpha_i)<\infty$, and Lemma \ref{lm:retraction*} gives $i(l(E_n),\alpha_i)\leq L$, so that the first claim follows from Lemma \ref{lm:compactML}.

If $c(E_n)\to c$, then $\mathrm{gr}_{l(E_n)}(m(E_n))=c(E_n)$ stays in a compact subset of $\cT_S$, and Corollary \ref{cr:joint*} shows that $(l(E_n),m(E_n))_{n\in \N}$ subconverges to some $(l,m)$. Let $E$ be the hyperbolic end with $\sigma(E)=\mathrm{Gr}_l(m)$. Then $E_n\to E$ and, by the continuity of grafting, $c(E)=\lim_{n\to \infty} c(E_n)=c$.
\end{proof}

Note that, as for Proposition \ref{pr:endsconverge}, neither the convexity of $u^*_n(S)$ nor the hypothesis $K^*<1$ is used here.

\subsubsection{A priori estimates}
\label{sssc:dualest}

We first translate the curvature hypothesis. By Lemma \ref{lm:dual}, $K^*\det B=K$ and $\det B=1+K$, so that
\begin{equation}\label{eq:KKstar}
  K\;=\;\frac{K^*}{1-K^*}~.
\end{equation}
Hence a sequence $h^*_n\to h^*$ with $K_{h^*}<0$ has $\inf K_{h_n}\to\inf K_{h^*}/(1-\inf K_{h^*})>-1$ and $\sup K_{h_n}\to\sup K_{h^*}/(1-\sup K_{h^*})<0$, where $h_n$ is the induced metric of the surface with third fundamental form $h^*_n$. Similarly, since $da=\det(B^*)\,da^*=(1-K^*)\,da^*$,
\begin{equation}\label{eq:areadual}
  \area(h_n)\;=\;\int_S(1-K^*_n)\,da^*_n\;=\;\area(h^*_n)-2\pi\chi(S)~,
\end{equation}
which is bounded. The diameter of $h_n$, on the other hand, is not controlled by $h^*_n$.

\begin{proposition}\label{pr:height*}
Let $(h^*_n)_{n\in \N}$ be a sequence of metrics on $S$ converging smoothly to $h^*$ with $K_{h^*}<0$, and let $(E_n,u^*_n)_{n\in \N}$ be a sequence of elements of $\cE$ where $u^*_n:S\to E_n$ is a good embedding with third fundamental form $h^*_n$, with $(E_n)_{n\in \N}$ converging to $E$. Let $h_n$ be the induced metric on $u^*_n(S)$. Then there are constants $C_1,C_2,C_3>0$, independent of $n$, such that
$$ \forall x\in u^*_n(S), C_1\leq d(x,\partial_0E_n) \leq C_2~, $$
$$ \int_SH^*_n da^*_n = \int_SH_n da_n \leq C_3~. $$
\end{proposition}

\begin{proof}
By \eqref{eq:KKstar} we can fix $k_0\leq\inf_n\inf_SK_{h_n}$ and $k\geq\sup_n\sup_SK_{h_n}$, both in $(-1,0)$.

\emph{Lower bound.} By Theorem \ref{tm:labourie41}(i), $u^*_n(S)$ lies outside $S_{k_0}(E_n)$, so that $d(x,\partial_0E_n)\geq d(S_{k_0}(E_n),\partial_0E_n)$ for all $x\in u^*_n(S)$; the right-hand side converges to $d(S_{k_0}(E),\partial_0E)>0$ by Lemma \ref{lm:kuniform}(ii).

\emph{Upper bound.} Since $K_{h_n}\leq k<0$, Lemma \ref{lm:kuniform}(i) applied to $\Sigma=u^*_n(S)$ gives
$$ d(x,\partial_0E_n)\leq\operatorname{artanh}\sqrt{1+k}=:C_2~. $$

\emph{Total mean curvature.} The two integrals coincide by (D5). The bound on $\int_SH_n\,da_n$ is Theorem \ref{tm:labourie41}(iv), whose proof in \cite[4.1(iv)]{L4} uses only $C_1$, $C_2$, a lower bound on $K_{h_n}$, an upper bound on $\area(h_n)$, and the geometry of $E_n$ on the region $\{C_1\leq d(\cdot,\partial_0E_n)\leq C_2\}$. By Lemma \ref{lm:geomconv} the injectivity radius of $E_n$ is bounded below on this region uniformly in $n$, so $u^*_n(S)$ is covered by a bounded number of balls isometric to balls of $\HH^3$, and in each of them \cite[Lemme~5.4(iii)]{L1} bounds the total mean curvature of a convex surface in terms of the radius and of $\inf K_{h_n}$. The bound on $\area(h_n)$ follows from \eqref{eq:areadual}.
\end{proof}

\subsubsection{Compactness}
\label{sssc:compact*}

Labourie's compactness theorem (Theorem \ref{tm:labourieD}) cannot be applied directly to the hyperbolic surfaces $u^*_n(S)$, since their induced metrics $h_n$ are not known to converge, nor even to be bounded, without a bound on $B_n$. We use instead its Lorentzian analogue on the dual surfaces, whose induced metrics are the $h^*_n$. Here ``uniformly elliptic'' means that $\det B^*\geq\varepsilon>0$, that is, $K^*\leq 1-\varepsilon$; it is the uniform version of the local strong convexity appearing in (D3), and holds trivially when $K^*<0$.

We will use the following statement.

\begin{theorem}[{\cite[Thm 5.5]{these}}]\label{tm:sch55}
Let $(\phi_n)_{n\in \N}$ be a sequence of uniformly elliptic spacelike immersions $\phi_n:D\to dS^3$ of the disc, whose induced metrics converge in $C^\infty$ to a metric $g_\infty$, such that $\int_DH_n\,da_n$ is bounded and, for some $x_n\to x_\infty$, $j^1\phi_n(x_n)$ converges. Then a subsequence converges in $C^\infty$ on compact subsets to an isometric immersion.
\end{theorem}

Note that no $C^0$ convergence is assumed: in the Lorentzian setting it is a consequence of the bound on the mean curvature, see the proof of \cite[Thm 5.5]{these}.

\begin{proposition}\label{pr:endsimplyconv*}
Assume that $h^*_n\to h^*$ smoothly, with $K_{h^*}<0$. Let $(E_n,u^*_n)_{n\in \N}$ be a sequence of elements of $\cE$ where $u^*_n:S\to E_n$ is a good embedding with third fundamental form $h^*_n$, with $(E_n)_{n\in \N}$ converging to $E$. Then, after extracting a subsequence, $(E_n,u^*_n)_{n\in \N}$ converges smoothly to $(E,u^*)\in\cE$, where $u^*:S\to E$ is a good embedding with third fundamental form $h^*$.
\end{proposition}

\begin{proof}
  Let $(\phi_n,\rho_n)$ be the lift of $u^*_n$ to an equivariant immersion of $\tilde S$ in $\HH^3$ given by Lemma \ref{lm:equiv}, and let $f_n=N_n:\tilde S\to dS^3$ be its unit normal, which by (D3) is a $\rho_n$-equivariant spacelike immersion with induced metric $\tilde h^*_n$ and shape operator $B^*_n=B_n^{-1}$, uniformly elliptic since $\det B^*_n=1-K^*_n\geq 1$. Compose with isometries of $dS^3$ (equivalently of $\HH^3$) so that $j^1f_n(x_0)$ converges for a fixed $x_0\in\tilde S$.
  
  By Proposition \ref{pr:height*}, $\int_SH^*_n\,da^*_n\leq C_3$, so that for every disc $D\subset\tilde S$ meeting $m$ fundamental domains of $\pi_1S$, $\int_DH^*_n\,da^*_n\leq mC_3$. Theorem \ref{tm:sch55} applied on a disc centered at $x_0$ gives a subsequence converging in $C^\infty$ on a smaller disc $D_0$; for $x_1\in D_0$, the $1$-jets $j^1f_n(x_1)$ then converge, and Theorem \ref{tm:sch55} applies again on a disc centered at $x_1$.

  Since $\tilde S$ is connected, every compact subset of $\tilde S$ is covered by finitely many discs reached from $x_0$ in this way by finite chains. Taking an exhaustion of $\tilde S$ by compact subsets and a diagonal extraction, we obtain $f_n\to f_\infty$ in $C^\infty_{\rm loc}(\tilde S)$, where $f_\infty$ is an isometric immersion of $(\tilde S,\tilde h^*)$ in $dS^3$, uniformly elliptic since $\det B^*=1-K^*\geq 1$.

  For each $\gamma\in\pi_1S$, $\rho_n(\gamma)$ is the isometry of $dS^3$ sending the $1$-jet of $f_n$ at $x_0$ to its $1$-jet at $\gamma x_0$, and since both $1$-jets converge, $\rho_n(\gamma)$ converges to an isometry $\rho'(\gamma)$, which defines a representation $\rho':\pi_1S\to\PSL(2,\C)$ for which $f_\infty$ is equivariant.

  By (D3), the unit normal $\phi$ of $f_\infty$ is an equivariant locally strongly convex immersion of $\tilde S$ in $\HH^3$ with third fundamental form $\tilde h^*$, and $\phi_n\to\phi$ in $C^\infty_{\rm loc}$ since $\phi_n$ is the unit normal of $f_n$. By Lemma \ref{lm:equiv}, $\phi$ defines a good embedding $u^*:S\to E'$ in a hyperbolic end $E'$, with third fundamental form $h^*$. The $\bCP^1$-structure of $E'$ is that of $(\phi,\rho')$, that is, the limit of the $\bCP^1$-structures of $(\phi_n,\rho_n)$, which are the $\sigma(E_n)\to\sigma(E)$, and Theorem \ref{tm:thurston} shows that $E'=E$.
\end{proof}

\subsubsection{Proof of Lemma \ref{lm:proper*}}
\label{sssc:proof*}

Let $h^*_n\to h^*$, $(E_n,u^*_n)_{n\in \N}$ and $c_n\to c$ be as in Lemma \ref{lm:proper*}. By Proposition \ref{pr:endsconverge*}, after extracting a subsequence, $E_n\to E$ with $c(E)=c$. By Proposition \ref{pr:endsimplyconv*}, after a further extraction, $(E_n,u^*_n)_{n\in \N}$ converges smoothly to $(E,u^*)$, where $u^*:S\to E$ is a good embedding with third fundamental form $h^*$. Since $c(E)=c$, this is the conclusion of Lemma \ref{lm:proper*}.
\qed

Note that the proof uses the hypothesis that $K^*<0$ only to obtain an upper bound on the restriction to $u^*(S)$ of the cosmological time of $E^*$ (or, dually, the distance to $\partial_0E$ in $E$). This is a shortcut, and it is quite conceivable that, with more effort, the hypothesis can be weakened to $K^*\leq 0$. In fact it might be possible (using in a more extensive way the compactness result in \cite{these}) to prove Lemma \ref{lm:proper*} under the hypothesis that $h^*$ has curvature $K^*<1$, with all closed, contractible geodesics of length $L>2\pi$. We do not attempt this here since the infinitesimal rigidity statement that we need (Lemma \ref{lm:rigid*}) requires $K^*\leq 0$ anyway.

\section{Proofs of the main results}

Throughout, $\cM_S$ is the space of smooth metrics on $S$ with curvature in $(-1,0]$, and $\cM^*_S$ the space of smooth metrics with curvature in $(-\infty,0)$, both with the $C^\infty$ topology. We implement here the sketch of proof in Section \ref{ssc:sketch}.

\subsection{The deformation spaces}\label{ssc:defspaces}

Recall that for $h\in \bar \cM_S$, $\cG_h$ is the space of $h$-self-adjoint positive definite $B\in\Gamma(\mathrm{End}\,TS)$ with $\det B=1+K_h$ and $d^DB=0$ (solutions of the Gauss--Codazzi equations), and $\cG=\{(h,B)~|~h\in\bar \cM_S, B\in\cG_h\}$. For $h^*\in\bar \cM^*_S$, $\cG^*_{h^*}$ is the space of $h^*$-self-adjoint $B^*$ with $\det B^*=1-K_{h^*}$ and $d^{D^*}B^*=0$, and $\cG^*=\{(h^*,B^*)~|~h^*\in\bar \cM^*_S, B^*\in\cG^*_{h^*}\}$.

By the fundamental theorem of surface theory and Lemma~\ref{lm:equiv}, both $\cG$ and $\cG^*$ are in natural bijection with $\cE$: $(h,B)$ maps to the pair $(E,u)$ whose induced metric is $h$ and shape operator $B$, and $(h^*,B^*)$ maps to the pair with third fundamental form $h^*$ and $B=(B^*)^{-1}$ (Lemma~\ref{lm:dual}). We give $\cE$ the topology of $\cG$ (equivalently of $\cG^*$: the bijection $(h,B)\mapsto(h(B\cdot,B\cdot),B^{-1})$ is a homeomorphism). ``Smooth convergence'' in Lemmas \ref{lm:proper} and \ref{lm:proper*} is convergence in this topology.

The conformal structure at infinity is an explicit function of these data,
$$ c(h,B)=[\,h((E+B)\cdot,(E+B)\cdot)\,]=[\,h^*((E+B^*)\cdot,(E+B^*)\cdot)\,]\in\cT_S~, $$
by Definition~\ref{df:Isharp} and the remark following it. The map sending $(h,B)$ to $c(h,B)$ is a smooth map because the projection from metrics to $\cT_S$ is smooth. We write
$$
\begin{array}{cccc}
  \Phi: & \cG & \to & \bar \cM_S\times\cT_S \\
        & (h,B) & \mapsto & (h,c(h,B))~, \\
  \Phi^*: & \cG^* & \to & \bar \cM^*_S\times\cT_S \\
  & (h^*,B^*)  & \mapsto & (h^*,c(h^*,B^*))~, 
\end{array} $$
and $\Phi_h=c|_{\cG_h}$, $\Phi^*_{h^*}=c|_{\cG^*_{h^*}}$. Theorem~\ref{tm:main} is the statement that $\Phi_h$ is bijective for every $h\in\cM_S$, and Theorem~\ref{tm:main*} that $\Phi^*_{h^*}$ is bijective for every $h^*\in\cM^*_S$.

\begin{theorem}[Labourie {\cite[Lemmes 3.1, 3.2]{L4}}; Labourie--Schlenker {\cite[Lemme 6.1]{iie}}]
  \label{tm:manifold}
\begin{enumerate}
\item For $h\in\cM_S$, $\cG_h$ is either empty or a smooth manifold of dimension
$6g-6$, whose tangent space at $B$ is
$\{\dot B\ h\text{-self-adjoint}:\ d^D\dot B=0,\ \tr(B^{-1}\dot B)=0\}$.
\item For a smooth path $(h_t)_{t\in[0,1]}$ in $\cM_S$,
$W_{(h_t)}=\{(B,t):B\in\cG_{h_t}\}$ is either empty or a smooth manifold with
boundary of dimension $6g-5$, with $\partial W_{(h_t)_{t\in[0,1]}}=\cG_{h_0}\times\{0\}\sqcup\cG_{h_1}\times\{1\}$.
\item The same holds for $\cG^*_{h^*}$, $h^*\in\cM^*_S$, and for $W^*_{(h^*_t)_{t\in[0,1]}}$.
\end{enumerate}
\end{theorem}

Points (1) and (2) are stated in \cite{L4} for $K>-1$, while (3) is \cite[Lemme 6.1]{iie} for spacelike immersions in $dS^3$ ($K^*<1$), where the solution space of the Gauss equation $\det B^*=1-K^*$ is shown to be a manifold on which the Codazzi operator is a submersion of corank $6g-6$. The path version follows as in \cite[Lemme 3.2]{L4}. In all cases the manifold structure is the one induced from the Fr\'echet space of sections, so $c$ is smooth on these manifolds.

\subsection{The covering map property}

We state here a form of the local diffeomorphism of the map $\Phi_h$ adapted to the argument below, which varies the metric.

\begin{proposition}\label{pr:localhomeo}
For $h\in\cM_S$, $\Phi_h\colon\cG_h\to\cT_S$ is a local diffeomorphism. For a smooth path $(h_t)_{t\in [0,1]}$ in $\cM_S$,
$$
\begin{array}{cccc}
  \Phi_{(h_t)_{t\in [0,1]}}: & W_{(h_t)_{t\in [0,1]}} & \to& \cT_S\times[0,1]~, \\
  & (B,t) & \mapsto & (c(h_t,B),t) 
\end{array}
$$
is a local diffeomorphism of manifolds with boundary. The same holds for $\Phi^*_{h^*}$ and $\Phi^*_{(h^*_t)_{t\in [0,1]}}$, $h^*_t\in\cM^*_S$.
\end{proposition}

\begin{proof}
Let $\dot B\in T_B\cG_h$ with $d\Phi_h(\dot B)=0$. The variation of $I^\#$ along $\dot B$ is then tangent to the orbit of $\mathrm{Diff}_0(S)\times C^\infty(S)$ acting by pull-back and conformal change, i.e.\ $\dot I^\#=\cL_vI^\#+2\lambda I^\#$ for some $v,\lambda$; this is precisely a first-order deformation of $(E,\Sigma)$ preserving the induced metric and the isotopy class of the conformal structure at infinity, in the sense of Section \ref{ssc:variations}. Since $K_h\le0$, Lemma~\ref{lm:rigid} gives $\dot B=0$. Thus $d\Phi_h$ is injective between spaces of dimension $6g-6$, hence an isomorphism, and the inverse function theorem applies. For $W_{(h_t)_{t\in [0,1]}}$: if $d\Phi_{(h_t)_{t\in [0,1]}}(\dot B,\dot t)=0$ then $\dot t=0$ and $\dot B\in T_B\cG_{h_t}$ with $d\Phi_{h_t}(\dot B)=0$, so $\dot B=0$, while dimensions are $6g-5$ on both sides. The dual statements use Lemma~\ref{lm:rigid*} and Theorem~\ref{tm:manifold}(3).
\end{proof}

The following  well-known statement will be used below.

\begin{lemma}\label{lm:covering}
Let $f\colon X\to Y$ be a proper local homeomorphism between locally compact
Hausdorff spaces, with $Y$ connected. Then $f$ is a covering map with finite
fibres, and $\deg f=\#f^{-1}(y)$ is independent of $y\in Y$.
\end{lemma}

This is standard, see e.g.\ \cite[Thm.~4.22]{forster}, where it is stated for locally compact spaces.

\begin{proposition}\label{pr:covering}
Let $(h_t)_{t\in[0,1]}$ be a smooth path in $\cM_S$. Then $\Phi_{(h_t)_{t\in [0,1]}}$ is proper, hence a finite covering of $\cT_S\times[0,1]$, and
$$ \deg\Phi_{h_0}=\deg\Phi_{(h_t)_{t\in [0,1]}}=\deg\Phi_{h_1}~. $$
As a consequence, $\deg\Phi_h$ is locally constant on $\cM_S$. The same holds for $\Phi^*$ over paths in $\cM^*_S$.
\end{proposition}

\begin{proof}
Let $(B_n,t_n)\in W_{(h_t)_{t\in[0,1]}}$ with $(c_n,t_n)=\Phi_{(h_t)_{t\in[0,1]}}(B_n,t_n)\to(c,t)$. Then $h_{t_n}\to h_t$ smoothly, so Lemma~\ref{lm:proper} applied to the pairs $(E_n,u_n)$ associated with $(h_{t_n},B_n)$ gives a subsequence converging in $\cE$ to $(E,u)$ with induced metric $h_t$ and $c(E)=c$, i.e. $(B_n,t_n)$ subconverges in $W_{(h_t)_{t\in[0,1]}}$. Since $\cT_S\times[0,1]$ is metrisable, this is properness. Lemma \ref{lm:covering} then applies (Proposition~\ref{pr:localhomeo}), and the fibre over $(c,i)$ is $\Phi_{h_i}^{-1}(c)$ for $i=0,1$. For $\Phi^*$ we use Lemma \ref{lm:proper*}, whose hypothesis $K^*<0$ is compatible with the definition of $\cM^*_S$.
\end{proof}

Note also that Lemma~\ref{lm:proper} does not require $K\le0$: it is the rigidity, not the properness, that uses the sign of the curvature.

\subsection{Connectedness of the spaces of metrics}

The argument also requires the fact that the spaces of metrics considered are connected.

\begin{lemma}\label{lm:connected}
$\cM_S$ and $\cM^*_S$ are connected. More precisely, every $h\in\cM_S$ can be
joined within $\cM_S$ to a metric of constant curvature in $(-1,0)$, and every
$h^*\in\cM^*_S$ can be joined within $\cM^*_S$ to a metric of constant curvature in
$(-\infty,0)$.
\end{lemma}

\begin{proof}
Write $h=e^{2u}g_c$ with $g_c$ hyperbolic in the conformal class of $h$, and let $s\in [0,1]$. The curvature of $h_s=e^{2su}g_c$ is
$$ K_s=e^{-2su}(-1-s\Delta_cu)=-(1-s)e^{-2su}+s\,K_h\,e^{2(1-s)u}~, $$
using $\Delta_cu=-1-K_he^{2u}$. Both terms are non-positive, so $K_s\le0$ for all $s\in[0,1]$, and $K_s\ge-C$ with $C=(1+\|K_h\|_\infty)e^{2\|u\|_\infty}$.

The path $(Ch_s)_{s\in [0,1]}$ stays in the space of metrics with curvature in $(-1,0]$ and connects $Cg_c$ (constant curvature $-1/C\in(-1,0)$) to $Ch$, while  the path $(\mu^2h)_{\mu\in[1,\sqrt{C}]}$ stays in the same space of metrics of curvature in $(-1,0]$ and connects $Ch$ to $h$. Constant curvature metrics with curvature in $(-1,0)$ form a connected set. This proves the claim for $\cM_S$.

For $\cM^*_S$, the path $(h_s)_{s\in [0,1]}$ itself stays in the space of metrics of negative curvature, since the first term of $K_s$ is negative for $s<1$ and $K_1=K_h<0$, and it ends at $g_c$, which can be rescaled to any negative constant curvature.
\end{proof}

\subsection{The degree is one}

We now consider specifically metrics of constant curvature.

\begin{theorem}[Labourie \cite{L5}, see {\cite[\S 6.1]{L4}}]\label{tm:fibreRn}
For $k\in(-1,0)$ and $g_k$ a metric of constant curvature $k$ on $S$, $\cG_{g_k}$
is homeomorphic to $\R^{6g-6}$; in particular it is nonempty and connected.
\end{theorem}

In Labourie's notation, $\cG_{g_k}=I(S,g_k)$ is the fibre over $[g_k]$ of the
projection $\pi_k\colon\cP_S\to\cT_S$ sending a $\bCP^1$-structure to the
conformal class of the induced metric on the $k$-surface of the corresponding
end; $\pi_k$ is a fibration with fibre homeomorphic to $\R^{6g-6}$.

\begin{proposition}\label{pr:degone}
For $g_k$ of constant curvature $k\in(-1,0)$, $\Phi_{g_k}\colon\cG_{g_k}\to\cT_S$ is
a homeomorphism. Consequently $\deg\Phi_h=1$ for every $h\in\cM_S$.
\end{proposition}

\begin{proof}
By Propositions~\ref{pr:localhomeo} and \ref{pr:covering},
$\Phi_{g_k}$ is a finite covering of $\cT_S$. Since $\cT_S$ is simply connected and
$\cG_{g_k}$ is connected and nonempty (Theorem~\ref{tm:fibreRn}), the covering is
trivial with one sheet, i.e.\ $\Phi_{g_k}$ is a homeomorphism and
$\deg\Phi_{g_k}=1$. For general $h\in\cM_S$, join $h$ to some $g_k$ by a smooth
path in $\cM_S$ (Lemma~\ref{lm:connected}) and apply Proposition~\ref{pr:covering}.
\end{proof}

For the third fundamental form the analogous input is:

\begin{theorem}[Labourie \cite{L5}, Schoen \cite{schoen:role}]\label{tm:minlag}
Let $g,g'$ be hyperbolic metrics on $S$. There is a unique $g$-self-adjoint, positive definite bundle morphism $b:TS\to TS$ with $\det b=1$ and $d^Db=0$ ($D$ the Levi-Civita connection of $g$) such that $g(b\cdot,b\cdot)$ is isotopic to $g'$.
\end{theorem}

Equivalently, there is a unique minimal Lagrangian diffeomorphism $(S,g)\to(S,g')$ isotopic to the identity: a diffeomorphism $m$ is minimal Lagrangian if and only if $m^*g'=g(b\cdot,b\cdot)$ with $b$ as above.

\begin{proposition}\label{pr:fibreRn*}
Let $k^*<0$ and let $g^*$ be a metric of constant curvature $k^*$ on $S$. Then $\cG^*_{g^*}$ is nonempty and connected. Therefore, $\Phi^*_{g^*}$ is a homeomorphism from $\cG^*_{g^*}$ to $\cT_S$. 
\end{proposition}

\begin{proof}
Let $B^*\in\cG^*_{g^*}$, so $B^*$ is $g^*$-self-adjoint, positive definite,
Codazzi for $D^*$, and $\det B^*=1-k^*$ is constant. Set $b^*=(1-k^*)^{-1/2}B^*$ and $\bar g=-k^*g^*$, which is hyperbolic and has the same Levi-Civita connection as $g^*$. Then $b^*$ is $\bar g$-self-adjoint, positive, Codazzi, with $\det b^*=1$. Conversely every such $b^*$ gives an element of $\cG^*_{g^*}$.

Theorem \ref{tm:minlag} with $g=\bar g$ shows that
$$ \begin{array}{cccc}
     \Psi: & \cG^*_{g^*}& \to & \cT_S~, \\
          & B^* & \mapsto & [\bar g(b^*\cdot,b^*\cdot)]
   \end{array} $$
is bijective. It is continuous, and $\cG^*_{g^*}$ is a manifold of dimension $6g-6$ by Theorem \ref{tm:manifold}(3), so $\Psi$ is a homeomorphism by invariance of domain. Hence $\cG^*_{g^*}\cong\R^{6g-6}$. 

Proposition \ref{pr:covering} then shows that the covering $\Phi^*_{g^*}$ of the simply connected space $\cT_S$ has one sheet.
\end{proof}

\begin{proposition}\label{pr:degone*}
$\deg\Phi^*_{h^*}=1$ for every $h^*\in\cM^*_S$.
\end{proposition}

\begin{proof}
Identical to Proposition~\ref{pr:degone}, with Lemma~\ref{lm:connected} for
$\cM^*_S$.
\end{proof}

\subsection{Proof of the main results}

\begin{proof}[Proof of Theorem \ref{tm:main} and Corollary \ref{cr:main}]
Let $h\in\cM_S$ and $c\in\cT_S$. By Propositions~\ref{pr:localhomeo}, \ref{pr:covering} and \ref{pr:degone}, $\Phi_h\colon\cG_h\to\cT_S$ is a covering of degree one, i.e.\ a homeomorphism. Hence there is exactly one $B\in\cG_h$ with $c(h,B)=c$, i.e.\ exactly one pair $(E,u)\in\cE$ with induced metric $h$ and conformal structure at infinity $c$. This is Corollary~\ref{cr:main}, and Theorem~\ref{tm:main} follows by Lemma~\ref{lm:equiv}.
\end{proof}

\begin{proof}[Proof of Theorem \ref{tm:main*} and Corollary \ref{cr:main*}]
  The argument is the same, with $\Phi^*_{h^*}$, Lemma~\ref{lm:rigid*}, Lemma~\ref{lm:proper*}, Theorem~\ref{tm:manifold}(3) and Proposition~\ref{pr:degone*}.
\end{proof}

\bibliographystyle{alpha}
\bibliography{/home/jean-marc/Dropbox/papiers/outils/biblio}

\def\cprime{$'$}
\begin{thebibliography}{KMT06b}

\bibitem[Ale05]{alex}
Alexander~D. Alexandrov.
\newblock {\em Convex polyhedra}.
\newblock Springer Monographs in Mathematics. Springer-Verlag, Berlin, 2005.
\newblock Translated from the 1950 Russian edition by N. S. Dairbekov, S. S.
  Kutateladze and A. B. Sossinsky, With comments and bibliography by V. A.
  Zalgaller and appendices by L. A. Shor and Yu. A. Volkov.

\bibitem[BB09]{benedetti-bonsante}
Riccardo Benedetti and Francesco Bonsante.
\newblock Canonical {Wick} rotations in 3-dimensional gravity.
\newblock {\em Memoirs of the American Mathematical Society}, 198:164pp, 2009.
\newblock math.DG/0508485.

\bibitem[BBZ11]{BBZ2}
Thierry Barbot, Fran{\c{c}}ois B{\'e}guin, and Abdelghani Zeghib.
\newblock Prescribing {Gauss} curvature of surfaces in 3-dimensional
  spacetimes, application to the {Minkowski} problem in {Minkowski} space.
\newblock {\em Ann. Inst. Fourier (Grenoble)}, 61(1):511--591, 2011.

\bibitem[Bel17]{belraouti2}
Mehdi Belraouti.
\newblock Asymptotic behavior of {C}auchy hypersurfaces in constant curvature
  space-times.
\newblock {\em Geom. Dedicata}, 190:103--133, 2017.

\bibitem[BMS13]{cyclic}
Francesco Bonsante, Gabriele Mondello, and Jean-Marc Schlenker.
\newblock A cyclic extension of the earthquake flow {I}.
\newblock {\em Geom. Topol.}, 17(1):157--234, 2013.

\bibitem[BMS15]{cyclic2}
Francesco Bonsante, Gabriele Mondello, and Jean-Marc Schlenker.
\newblock A cyclic extension of the earthquake flow {II}.
\newblock {\em Ann. Sci. \'Ec. Norm. Sup\'er. (4)}, 48(4):811--859, 2015.

\bibitem[DW08]{dumas-wolf}
Emily Dumas and Michael Wolf.
\newblock Projective structures, grafting and measured laminations.
\newblock {\em Geom. Topol.}, 12(1):351--386, 2008.

\bibitem[Eps24]{epstein:envelopes}
Charles~L Epstein.
\newblock Envelopes of horospheres and {Weingarten} surfaces in hyperbolic
  3-space.
\newblock {\em arxiv:2401.12115}, 2024.
\newblock Original preprint from 1984.

\bibitem[Fil07]{fillastre2}
Fran{\c{c}}ois Fillastre.
\newblock Polyhedral realisation of hyperbolic metrics with conical
  singularities on compact surfaces.
\newblock {\em Ann. Inst. Fourier (Grenoble)}, 57(1):163--195, 2007.

\bibitem[For81]{forster}
Otto Forster.
\newblock {\em Lectures on {R}iemann Surfaces}, volume~81 of {\em Graduate
  Texts in Mathematics}.
\newblock Springer-Verlag, New York, 1981.
\newblock Translated from the German by Bruce Gilligan.

\bibitem[FS16]{fillastre-seppi}
Fran{\c{c}}ois Fillastre and Andrea Seppi.
\newblock Spherical, hyperbolic and other projective geometries: convexity,
  duality, transitions.
\newblock {\em arXiv preprint arXiv:1611.01065}, 2016.

\bibitem[FS19]{fillastre-seppi:spherical}
Fran\c{c}ois Fillastre and Andrea Seppi.
\newblock Spherical, hyperbolic, and other projective geometries: convexity,
  duality, transitions.
\newblock In {\em Eighteen essays in non-{E}uclidean geometry}, volume~29 of
  {\em IRMA Lect. Math. Theor. Phys.}, pages 321--409. Eur. Math. Soc.,
  Z\"{u}rich, 2019.

\bibitem[Gro86]{PDR}
M.~Gromov.
\newblock {\em Partial Differential Relations}.
\newblock Springer, 1986.

\bibitem[HR93]{HR}
Craig~D. Hodgson and Igor Rivin.
\newblock A characterization of compact convex polyhedra in hyperbolic 3-space.
\newblock {\em Invent. Math.}, 111:77--111, 1993.

\bibitem[KMT03]{KMT}
Sadayoshi Kojima, Shigeru Mizushima, and Ser~Peow Tan.
\newblock Circle packings on surfaces with projective structures.
\newblock {\em J. Differential Geom.}, 63(3):349--397, 2003.

\bibitem[KMT06a]{KMT3}
Sadayoshi Kojima, Shigeru Mizushima, and Ser~Peow Tan.
\newblock Circle packings on surfaces with projective structures: a survey.
\newblock In {\em Spaces of {K}leinian groups}, volume 329 of {\em London Math.
  Soc. Lecture Note Ser.}, pages 337--353. Cambridge Univ. Press, Cambridge,
  2006.

\bibitem[KMT06b]{KMT2}
Sadayoshi Kojima, Shigeru Mizushima, and Ser~Peow Tan.
\newblock Circle packings on surfaces with projective structures and
  uniformization.
\newblock {\em Pacific J. Math.}, 225(2):287--300, 2006.

\bibitem[KP94]{kulkarni-pinkall}
Ravi~S. Kulkarni and Ulrich Pinkall.
\newblock A canonical metric for {M}\"obius structures and its applications.
\newblock {\em Math. Z.}, 216(1):89--129, 1994.

\bibitem[KS08]{volume}
Kirill Krasnov and Jean-Marc Schlenker.
\newblock On the renormalized volume of hyperbolic 3-manifolds.
\newblock {\em Comm. Math. Phys.}, 279(3):637--668, 2008.

\bibitem[KS09]{cp}
Kirill Krasnov and Jean-Marc Schlenker.
\newblock A symplectic map between hyperbolic and complex {T}eichm\"uller
  theory.
\newblock {\em Duke Math. J.}, 150(2):331--356, 2009.

\bibitem[KT92]{kamishima-tan}
Yoshinobu Kamishima and Ser~P. Tan.
\newblock Deformation spaces on geometric structures.
\newblock In {\em Aspects of low-dimensional manifolds}, volume~20 of {\em Adv.
  Stud. Pure Math.}, pages 263--299. Kinokuniya, Tokyo, 1992.

\bibitem[Lab89]{L1}
Fran\c{c}ois Labourie.
\newblock Immersions isom\'etriques elliptiques et courbes
  pseudo-holo\-morphes.
\newblock {\em J. Differential Geom.}, 30:395--424, 1989.

\bibitem[Lab91]{L6}
Fran{\c{c}}ois Labourie.
\newblock Probl\`eme de {M}inkowski et surfaces \`a courbure constante dans les
  vari\'et\'es hyperboliques.
\newblock {\em Bull. Soc. Math. France}, 119(3):307--325, 1991.

\bibitem[Lab92a]{L4}
Fran\c{c}ois Labourie.
\newblock M\'etriques prescrites sur le bord des vari\'et\'es hyperboliques de
  dimension 3.
\newblock {\em J. Differential Geom.}, 35:609--626, 1992.

\bibitem[Lab92b]{L5}
Fran\c{c}ois Labourie.
\newblock Surfaces convexes dans l'espace hyperbolique et {CP1}-structures.
\newblock {\em J. London Math. Soc., II. Ser.}, 45:549--565, 1992.

\bibitem[LS00]{iie}
Fran{\c{c}}ois Labourie and Jean-Marc Schlenker.
\newblock Surfaces convexes fuchsiennes dans les espaces lorentziens \`a
  courbure constante.
\newblock {\em Math. Ann.}, 316(3):465--483, 2000.

\bibitem[Mes07]{mess}
Geoffrey Mess.
\newblock Lorentz spacetimes of constant curvature.
\newblock {\em Geom. Dedicata}, 126:3--45, 2007.

\bibitem[MS84]{morgan-shalen}
John~W. Morgan and Peter~B. Shalen.
\newblock Valuations, trees, and degenerations of hyperbolic structures. {I}.
\newblock {\em Ann. Math. (2)}, 120:401--476, 1984.

\bibitem[MST26]{grothu}
Daniel Monclair, Jean-Marc Schlenker, and Nicolas Tholozan.
\newblock Gromov-{Thurston} manifolds and anti-de {Sitter} geometry.
\newblock {\em Geom. Topol.}, 30(6):2097--2156, 2026.

\bibitem[Pog73]{Po}
Aleksei~V. Pogorelov.
\newblock {\em Extrinsic Geometry of Convex Surfaces}.
\newblock American Mathematical Society, 1973.
\newblock Translations of Mathematical Monographs. Vol. 35.

\bibitem[Sca99]{scannell}
Kevin~P. Scannell.
\newblock Flat conformal structures and the classification of de {S}itter
  manifolds.
\newblock {\em Comm. Anal. Geom.}, 7(2):325--345, 1999.

\bibitem[Sch93]{schoen:role}
Richard~M. Schoen.
\newblock The role of harmonic mappings in rigidity and deformation problems.
\newblock In {\em Complex geometry ({O}saka, 1990)}, volume 143 of {\em Lecture
  Notes in Pure and Appl. Math.}, pages 179--200. Dekker, New York, 1993.

\bibitem[Sch96]{these}
Jean-Marc Schlenker.
\newblock Surfaces convexes dans des espaces lorentziens \`a courbure
  constante.
\newblock {\em Comm. Anal. Geom.}, 4(1-2):285--331, 1996.

\bibitem[Sch98]{shu}
Jean-Marc Schlenker.
\newblock M\'etriques sur les poly\`edres hyperboliques convexes.
\newblock {\em J. Differential Geom.}, 48(2):323--405, 1998.

\bibitem[Sch02]{horo}
Jean-Marc Schlenker.
\newblock Hypersurfaces in {$H\sp n$} and the space of its horospheres.
\newblock {\em Geom. Funct. Anal.}, 12(2):395--435, 2002.

\bibitem[Sch06]{hmcb}
Jean-Marc Schlenker.
\newblock Hyperbolic manifolds with convex boundary.
\newblock {\em Invent. Math.}, 163(1):109--169, 2006.

\bibitem[Sch20]{weylsurvey}
Jean-Marc Schlenker.
\newblock On the {Weyl} problem for complete surfaces in the hyperbolic and
  anti-de {Sitter} spaces.
\newblock 2020.
\newblock arxiv:2012.14754. To appear, {\em Proceedings of the 14th Seasonal
  Institute of the Mathematical Society of Japan. Advanced Studies in Pure
  Mathematics.}

\bibitem[SW02]{scannell-wolf}
Kevin~P. Scannell and Michael Wolf.
\newblock The grafting map of {T}eichm\"uller space.
\newblock {\em J. Amer. Math. Soc.}, 15(4):893--927 (electronic), 2002.

\bibitem[SY18]{delaunay}
Jean-Marc. {Schlenker} and Andrew {Yarmola}.
\newblock {Properness for circle packings and Delaunay circle patterns on
  complex projective structures}.
\newblock {\em ArXiv e-prints}, June 2018.

\end{thebibliography}
\end{document}